\documentclass[a4paper,12pt,reqno]{amsart}
\usepackage[T1]{fontenc}
\usepackage[utf8]{inputenc}
\usepackage[english]{babel}
\usepackage[left=25mm,right=25mm,top=35mm,bottom=35mm]{geometry}
\usepackage{amsmath,amssymb,amsthm}
\usepackage{pgfplots}
\usepackage{pgfplotstable}
\usepackage[margin=0pt,font={small}]{caption}
\usepackage{booktabs}
\usepackage{color,colortbl}
\usepackage[hidelinks]{hyperref}

\newcommand{\eps}{\varepsilon}
\newcommand{\jj}{\mathbf{j}} 
\newcommand{\y}{\mathbf{y}} 
\newcommand{\z}{\mathbf{z}} 
\newcommand{\g}{\gamma}
\newcommand{\G}{\Gamma}
\newcommand{\N}{\mathbb{N}}
\renewcommand{\P}{\mathbb{P}}
\newcommand{\R}{\mathbb{R}}
\newcommand{\V}{\mathbb{V}}
\newcommand{\X}{\mathbb{X}}
\newcommand{\Y}{\mathbb{Y}}

\newcommand{\MM}{\mathcal{M}}
\newcommand{\NN}{\mathcal{N}}

\renewcommand{\SS}{\mathcal{S}}
\newcommand{\TT}{\mathcal{T}}

\DeclareMathOperator*{\hull}{span}
\DeclareMathOperator*{\refine}{refine}

\newcommand{\norm}[3][]{#1\|#2#1\|_{#3}}

\newcommand{\dx}{\mathrm{d}x}

\newcommand{\reff}[2]{\stackrel{\eqref{#1}}{#2}}	
\newcommand{\refp}[2]{\stackrel{\phantom{\eqref{#1}}}{#2}}
\newcommand{\FoneD}{\mathcal{F}} 
\newcommand{\Colpts}{\mathcal{Y}} 
\newcommand{\mf}{\kappa} 
\newcommand{\bmf}{\boldsymbol{\kappa}} 
\newcommand{\indset}{\Lambda} 
\newcommand{\markindset}{\Upsilon} 
\newcommand{\marg}{{\rm K}} 
\newcommand{\rmarg}{{\rm R}} 
\newcommand{\nnu}{\boldsymbol{\nu}} 
\newcommand{\eeps}{\boldsymbol{\eps}} 
\newcommand\1{\boldsymbol{1}}
\newcommand{\Lagr}{I} 
\newcommand{\LagrBasis}[2]{L_{#1}^{#2}} 
\newcommand{\LagrBasisHat}[2]{\widehat L_{#1}^{#2}} 
\newcommand{\scsol}{u_{\bullet}^{\rm SC}} 
\newcommand{\scsolhat}{\widehat u_{\bullet}^{\rm\, SC}} 

\pgfplotsset{
every axis/.append style={
font={\fontsize{8pt}{12pt}\selectfont},  
},
tick label style={font=\tiny},
title style={font=\tiny,yshift=-1.5ex},
xlabel style={font=\tiny,yshift=+1.0ex},
ylabel style={font=\tiny,yshift=-1.2ex},
}

\definecolor{myBrown}{rgb}{0.6 0.4 0.2}
\definecolor{myOrange}{rgb}{1.0 0.6 0.2}
\definecolor{myLightGray}{RGB}{235,235,235}
\definecolor{myViolet}{RGB}{153,50,204}
\newtheorem{theorem}{Theorem}

\newtheorem{algorithm}[theorem]{Algorithm}
\newtheorem{remark}[theorem]{Remark}

\makeatletter
\def\@seccntformat#1{%
  \protect\textup{\protect\@secnumfont
    \ifnum\pdfstrcmp{subsection}{#1}=0 \bfseries\fi
    \csname the#1\endcsname
    \protect\@secnumpunct
  }%
}
\makeatother
\newcommand*\patchAmsMathEnvironmentForLineno[1]{%
  \expandafter\let\csname old#1\expandafter\endcsname\csname #1\endcsname
  \expandafter\let\csname oldend#1\expandafter\endcsname\csname end#1\endcsname
  \renewenvironment{#1}%
     {\linenomath\csname old#1\endcsname}%
     {\csname oldend#1\endcsname\endlinenomath}}%
\newcommand*\patchBothAmsMathEnvironmentsForLineno[1]{%
  \patchAmsMathEnvironmentForLineno{#1}%
  \patchAmsMathEnvironmentForLineno{#1*}}%
\AtBeginDocument{%
\patchBothAmsMathEnvironmentsForLineno{equation}%
\patchBothAmsMathEnvironmentsForLineno{align}%
\patchBothAmsMathEnvironmentsForLineno{flalign}%
\patchBothAmsMathEnvironmentsForLineno{alignat}%
\patchBothAmsMathEnvironmentsForLineno{gather}%
\patchBothAmsMathEnvironmentsForLineno{multline}%
}
\usepackage[mathlines]{lineno}

\usepackage{fancyhdr}
\advance\footskip0.5cm
\newcommand\rev[1]{{\color{black}#1}}
\definecolor{otherblue}{rgb}{0,0.3,0.6}
\def\rbl#1{\textcolor{black}{#1}}
\def\rblx#1{\textcolor{otherblue}{#1}}
\newcommand\abrev[1]{{\color{black}#1}}

\newcommand\abrevx[1]{{\color{black}#1}}

\title{Error estimation and adaptivity for stochastic collocation finite elements\abrevx{\\ Part I: }
single-level approximation}
\author{Alex Bespalov}
\address{School of Mathematics, University of Birmingham, Edgbaston, Birmingham B15 2TT, UK}
\email{a.bespalov@bham.ac.uk}
\author{David J. Silvester}
\address{Department of Mathematics, University of Manchester, Oxford Road, Manchester M13 9PL, UK}
\email{d.silvester@manchester.ac.uk}
\author{Feng Xu}
\address{Department of Mathematics, University of Manchester, Oxford Road, Manchester M13 9PL, UK}
\email{feng.xu@manchester.ac.uk}
\thanks{{\em Acknowledgements.}
\rev{This work was supported by EPSRC grants EP/W010925/1, EP/P013317/1, and~EP/P013791/1.}
}
\date{\today}

\begin{document}

\begin{abstract}
A  general adaptive refinement strategy for solving linear elliptic partial differential equation with random data is proposed and analysed herein. The adaptive strategy \abrevx{extends} the a posteriori error estimation framework introduced by  Guignard \& Nobile in 2018 ({\sl SIAM J. Numer. Anal.}, {\bf 56}, 3121--3143) to cover  problems with  a {\it nonaffine} parametric coefficient dependence. A \abrevx{suboptimal}, but nonetheless reliable and convenient implementation of the strategy involves approximation of the decoupled  PDE problems with a common finite element approximation space.  Computational results obtained using  such a  {\it single-level} strategy are presented in this paper (part I). Results obtained using a potentially more efficient \abrevx{{\it multilevel}} approximation strategy, where meshes are individually tailored,  will be discussed in part~II of this work. The \rbl{codes} used to generate the numerical results \abrevx{are} available \rblx{on GitHub}.
\end{abstract}

\maketitle
\thispagestyle{fancy}

\section{Introduction} \label{sec:intro}

Partial differential equations (PDEs) with uncertain inputs feature prominently when modelling a host of physical phenomena 
and have become a de facto model over the last two decades, both in scientific computing and  computational engineering.
Sparse grid stochastic collocation representations of parametric uncertainty in combination with finite element discretization 
of physical space \rbl{have} become established as an alternative approach to Monte-Carlo strategies over the last decade, 
especially in the context of nonlinear PDE models or linear PDE problems that are nonlinear in the parameterisation of 
the uncertainty.

Sparse grid methods, where the set of sample points is adaptively generated, can be traced back to 
Gerstner \& Griebel~\cite{GerstnerGriebel2003}. They have been extensively tested \rbl{in} a collocation setting; see for 
example,~\cite{ChkifaCohenSchwab2014,NobileTamelliniTeseiTempone2016}. A complementary  concept  that has shown a lot 
of promise  is the employment of multilevel approaches that aim to reduce the computational cost through a hierarchy of spatial 
approximations; see for example,~\cite{LangSS20,TeckentrupJantschWebsterGunzburger2012}.
In this contribution, \rbl{we aim to} combine these two complementary concepts in a rigorous manner with adaptivity driven 
by \rbl{novel} reliable a posteriori error estimates.  Specifically, we will generalise and extend  the  adaptive framework 
\rbl{proposed} in the recent paper by  Guignard \& Nobile~\cite{GuignardN18} and  present a critical comparison of alternative 
strategies in the context of solving a representative model problem that combines strong anisotropy in the parametric dependence 
with singular behaviour in the physical  space. We note that parametric adaptivity has also been explored in  a Galerkin 
framework; see \cite{bps14,br18,bs16,egsz14, egsz15}\abrevx{,} and that there are a number of recent papers aimed at proving 
dimension-independent convergence; see for example,~\cite{BachmayrCohenDungSchwab2017,bprr18,NobileTamelliniTempone2016,ZechDungSchwab2018}.

\rbl{The convergence} of a modified version of the adaptive algorithm in~\cite{GuignardN18} has been established by 
Eigel et al.~\cite{eest22} and independently by Feischl \& Scaglioni~\cite{FeischlS21}. The authors of~\cite{FeischlS21} note 
that the main obstacle in establishing convergence  is ``the interplay of parametric refinement and finite element refinement''. 
We focus on this interplay in this work. Specifically, after \abrevx{introducing the model problem in section~\ref{sec:problem}
and setting up its discretization in section~\ref{sec:scfem}}, \rbl{we develop a general adaptive solution strategy  in 
section~\ref{sec:error:estimate}. Computational results generated with \abrevx{adaptively refined \emph{single-level}
 approximations (i.e., using the same spatial refinement for all collocation points)}
are discussed in section~\ref{sec:results}.  A discussion of computational results obtained with a more efficient
\abrevx{\emph{multilevel}} spatial refinement implementation of the adaptive algorithm in~\S\ref{sec:adaptivity}  is deferred 
to part~II of this work.}

\section{\rbl{A parametric model problem}} \label{sec:problem}

Let \rev{$D \subset \R^2$} be a bounded Lipschitz domain with polygonal boundary $\partial D$.
Let $\Gamma := \G_1 \times \G_2 \ldots \times \G_M$ denote the parameter domain in $\R^M$,
where $M \in \N$ and each $\G_m$ ($m = 1,\ldots,M$) is a bounded interval in~$\R$.
We introduce a probability measure $\pi(\y) := \prod_{m=1}^M \pi_m(y_m)$ on $(\G,\mathcal{B}(\G))$;
here, $\pi_m$ denotes a Borel probability measure on $\G_m$ ($m = 1,\ldots,M$) and
$\mathcal{B}(\G)$ is the Borel $\sigma$-algebra on $\G$.

We consider the following parametric elliptic problem:
find $u : \overline D \times \G \to \R$ satisfying
\begin{equation} 
\label{eq:pde:strong}
\begin{aligned}
-\nabla \cdot (a(\cdot, \y)\nabla u(\cdot, \y))
&= f 
&& \text{in $D$},\\ 
u(\cdot, \y) &= 0  && \text{on $\partial D$} 
\end{aligned}
\end{equation} 
$\pi$-almost everywhere on $\G$ (i.e., almost surely).
Here, the deterministic right-hand side function $f \in L^2(D)$ and
the coefficient $a$ is a random field on $(\G,\mathcal{B}(\G),\pi)$ over $L^\infty(D)$.
Furthermore, we assume that there exist constants $a_{\min},\, a_{\max}$ such that
\begin{equation} \label{eq:amin:amax}
   0 < a_{\min} \le \operatorname*{ess\;inf}_{x \in D} a(x,\y) \le 
   \operatorname*{ess\;sup}_{x \in D} a(x,\y) \le a_{\max} < \infty \quad
   \text{$\pi$-a.e. on $\G$}.
\end{equation}
This assumption, in par\-ti\-cu\-lar, implies the following norm equivalence:
for any $v \in \X := H^1_0(D)$ there holds
\begin{equation} \label{eq:norm:equiv}
   a_{\min}^{1/2} \|\nabla v\|_{L^2(D)} \le
   \| a^{1/2}(\cdot,\y) \nabla v \|_{L^2(D)} \le
   a_{\max}^{1/2} \|\nabla v\|_{L^2(D)}\quad
   \text{$\pi$-a.e. on $\G$}.
\end{equation}

The parametric problem~\eqref{eq:pde:strong} is understood in the \rbl{weak} sense:
given $f \in L^2(D)$, find $u : \G \to \X$ such that
\begin{align} \label{eq:pde:weak}
   \int_D a(x, \y) \nabla u(x,\y) \cdot \nabla v(x) \dx = \int_D f(x) v(x) \dx 
   \quad \forall v \in \X,\ \text{$\pi$-a.e. on $\G$}.
\end{align}
The above assumptions on $a$ and $f$ guarantee that the parametric problem~\eqref{eq:pde:strong}
admits a unique weak solution $u$ in the Bochner space $\V := L_\pi^p(\G; \X)$ for any $p \in [1, \infty]$;
\rbl{see~\cite[Lemma~1.1]{BabuskaNT07} for details}.

\section{Multilevel stochastic collocation finite element method} \label{sec:scfem}

For the numerical solution of problem~\eqref{eq:pde:strong} we \rbl{propose to} use the multilevel 
stochastic collocation finite element method (SC-FEM).
\rbl{We recall the main ideas and the construction of the approximation spaces in the following.}

Let $\TT_\bullet$ be a mesh, i.e., a conforming triangulation of the spatial domain $D$ into
compact non-degenerate triangles $T \in \TT_\bullet$ and denote by $\NN_\bullet$ the set of vertices of $\TT_\bullet$.
We \rbl{restrict attention to the space of continuous piecewise linear finite elements for convenience,}
\begin{equation*}
 \X_\bullet := \SS^1_0(\TT_\bullet) :=
 \{v \in \X : v \vert_T \text{ is affine for all } T \in \TT_\bullet \} \subset \X = H^1_0(D).
\end{equation*}
Recall that the standard basis of $\X_\bullet$ is given by
$\{ \varphi_{\bullet,\xi} : \xi \in \NN_\bullet \setminus \partial D \}$, where
$\varphi_{\bullet,\xi}$ denotes the hat function associated with the vertex $\xi \in \NN_\bullet$.

For mesh refinement, we employ newest vertex bisection (NVB); see, e.g., \cite{stevenson,kpp}.
We assume that any mesh $\TT_\bullet$ employed for the spatial discretization
can be obtained by applying NVB refinement(s) to a given (coarse) initial mesh $\TT_0$.

For a given mesh $\TT_\bullet$, 
let $\widehat\TT_\bullet$ be the coarsest NVB refinement of $\TT_\bullet$ such that 
all edges of $\TT_\bullet$ have been bisected once
(which corresponds to uniform refinement of all elements by three bisections). 
Then, $\widehat\NN_\bullet$ denotes the set of vertices of $\widehat\TT_\bullet$,
and $\NN_\bullet^+ := (\widehat\NN_\bullet \setminus \NN_\bullet) \setminus \partial D$
is the set of new interior vertices created by this refinement of $\TT_\bullet$.
The finite element space associated with $\widehat\TT_\bullet$ is denoted by
$\widehat\X_\bullet := \SS^1_0(\widehat\TT_\bullet)$,
and
$\{ \widehat\varphi_{\bullet,\xi} : \xi \in \widehat\NN_\bullet \setminus \partial D \}$ is the corresponding
 basis of hat functions. In~\S\ref{sec:error:estimate}, we will exploit the ($H^1$-stable) decomposition
\begin{equation} \label{eq:Y:space}
   \widehat\X_\bullet = \X_\bullet \oplus \Y_\bullet, \text{ \ where \ }
   \Y_\bullet := \hull\{ \widehat\varphi_{\bullet,\xi} : \xi \in \NN_\bullet^+ \}.
\end{equation}
Note that $\X_\bullet \cap \Y_\bullet = \{0\}$, therefore
the strengthened Cauchy--Schwarz inequality holds for the subspaces $\X_{\bullet}$ and $\Y_{\bullet}$
(see, e.g.~\cite{eijkhoutVass1991}):
\begin{equation} \label{eq:sCS}
   \exists\, \g \in [0,1) \text{ such that }
   \big| ( \nabla u, \nabla v )_{L^2(D)} \big|\le
   \g\, \norm{\nabla u}{L^2(D)}\, \norm{\nabla v}{L^2(D)}\ \
   \forall\, u \in \X_{\bullet},\ \forall\, v \in \Y_{\bullet}.
\end{equation}
For a set of marked vertices $\MM_\bullet \subseteq \NN_\bullet^+$,
let $\TT_\circ := \refine(\TT_\bullet,\MM_\bullet)$ be the coarsest NVB refinement of $\TT_\bullet$ such that
$\MM_\bullet \subset \NN_\circ$, i.e., all marked vertices are vertices of $\TT_\circ$.

For a fixed $\z \in \G$, consider a mesh $\TT_{\bullet \z}$ and its uniform refinement $\widehat\TT_{\bullet \z}$
as well as the corresponding finite element spaces
$\X_{\bullet \z} := \SS^1_0(\TT_{\bullet \z})$ and $\widehat\X_{\bullet \z} := \SS^1_0(\widehat\TT_{\bullet \z})$.
We denote by $u_{\bullet \z} \in \X_{\bullet \z}$ the Galerkin finite element solution satisfying
\begin{align} \label{eq:sample:fem}
   \int_D a(x, \z) \nabla u_{\bullet \z}(x) \cdot \nabla v(x) \dx = \int_D f(x) v(x) \dx 
   \quad \forall v \in \X_{\bullet \z}.
\end{align}
The enhanced Galerkin solution satisfying~\eqref{eq:sample:fem} for all $v \in \widehat\X_{\bullet \z}$
is denoted by $\widehat u_{\bullet \z} \in \widehat\X_{\bullet \z}$. 

Turning now to the parameter domain $\Gamma$,
we consider a finite set $ \Colpts_\bullet$ of  collocation points in $\G$.
The SC-FEM approximation of the solution $u$ to parametric problem~\eqref{eq:pde:strong} is built as
\begin{equation} \label{eq:scfem:sol}
   \scsol := \sum\limits_{\z \in \Colpts_\bullet} u_{\bullet \z}(x) \LagrBasis{\bullet \z}{}(\y),
\end{equation}
where $u_{\bullet \z} \in \X_{\bullet \z}$ are Galerkin approximations satisfying~\eqref{eq:sample:fem} 
for $\z \in \Colpts_\bullet$, and
$\{ \LagrBasis{\bullet \z}{}(\y) = \LagrBasis{\z}{\Colpts_\bullet}(\y) : \z \in \Colpts_\bullet \}$ is a set of
multivariable Lagrange basis functions associated with $\Colpts_\bullet$ and satisfying
$\LagrBasis{\bullet \z}{}(\z') = \delta_{\z\z'}$ for any $\z,\,\z' \in \Colpts_\bullet$.
The total number of degrees of freedom in the SC-FEM approximations defined by~\eqref{eq:scfem:sol}
 is given by $\sum_{\z \in \Colpts_\bullet} \dim(\X_{\bullet \z})$.
Note that the SC-FEM solution considered here 
follows the so-called \emph{multilevel} construction
(cf.~\cite{LangSS20,FeischlS21}) that allows $\X_{\bullet \z} \not= \X_{\bullet \z'}$ for $\z \not= \z'$.
This is in contrast to the \emph{single-level} SC-FEM approximations that employ the same finite element space $\X_{\bullet}$
for all collocation points $\z \in \Colpts_\bullet$; see, e.g.,~\cite{BabuskaNT07,NobileTW08a,GuignardN18}.

Clearly, the choice of collocation points and the associated polynomial spaces on $\G$
is critical for efficient implementation of the generic SC-FEM construction outlined above,
particularly, for high-dimensional parametric problems.
The \rbl{established} methodology here utilizes the sparse grid idea that goes back to Smoljak in~\cite{Smoljak63}
\rbl{that is briefly described in the next \abrevx{section}}.

\subsection{Sparse grid interpolation} \label{sec:sparse:grids}

To simplify \rbl{the} presentation we assume that $\G_1 = \G_2 = \ldots = \G_M \subset \R$. 
The methodology extends trivially \rbl{to the general case}.
In order to construct a sparse grid $\Colpts_\bullet \subset \G = \G_1 \times \ldots \times \G_M$, one needs three ingredients:
\begin{itemize}
\item
a family $\FoneD$ of {\it nested sets} of 1D nodes on $\G_m$ (one family for all $m=1,\ldots,M$); 
examples of such node sets are Leja points and \rev{Clenshaw--Curtis quadrature points;}

\item
a strictly increasing function $\mf : \N_0 \to \N_0$ associated with the chosen sets of 1D nodes and such that
$\mf(0) = 0$, $\mf(1) = 1$ 
(e.g., $\mf(i) = i$ for Leja points and $\mf(i) = 2^{i-1}+1$, $i > 1$, for Clenshaw--Curtis nodes with the 
\rbl{usual} doubling rule);

\item
a monotone (or, downward-closed) finite set $\indset_\bullet \subset \N^M$ of multi-indices, i.e.,
$\indset_\bullet = \{ \nnu = (\nu_1,\ldots,\nu_M) : \nu_m \in \N\rbl{, \forall\, m = 1,\ldots,M} \}$ is such 
that $\#\indset_\bullet < \infty$ and
\[
   \nnu \in \indset_\bullet \Longrightarrow \nnu - \eeps_m \in \indset_\bullet \quad
   \forall\,m=1,\ldots,M \text{ such that } \nu_m >1,
\]
where $\eeps_m$ denotes the \rbl{$m$th} unit multi-index, i.e., $(\eeps_m)_i = \delta_{mi}$ 
for all $i=1,\ldots,M$.
Note that any monotone index set $\indset_\bullet$ contains the multi-index $\1 = (1, 1,\dots,1)$.
\end{itemize}

Now, for each $\nnu \in \indset_\bullet$, the set of collocation points along the  \rbl{$m$th} coordinate 
axis in $\R^M$ is given by \rev{the set $\Colpts_m^{\mf(\nu_m)} \in \FoneD$ of cardinality $\mf(\nu_m)$}
and we define
\[
   \Colpts^{\,(\nnu)} := \Colpts_1^{\mf(\nu_1)} \times \Colpts_2^{\mf(\nu_2)} 
   \times \ldots \times \Colpts_M^{\mf(\nu_M)}.
\]
For a given index set $\indset_\bullet$, the \emph{sparse grid} $\Colpts_\bullet = \Colpts_{\indset_\bullet}$
of collocation points on $\G$ is defined~as
\[
   \Colpts_\bullet = \Colpts_{\indset_\bullet} := \bigcup_{\nnu \in \indset_\bullet} \Colpts^{\,(\nnu)}.
\]

Let $\P_q$  denote the set of univariate polynomials of degree at most $q \in \N_0$.
Given an index set $\indset_\bullet$, we define the 
polynomial space $\P_{\bullet} = \P_{\indset_\bullet}$ on $\G$~as
\[
    \P_{\bullet} = \P_{\indset_\bullet} := \bigoplus_{\nnu \in \indset_\bullet} \P_{\bmf(\nnu)-\1}
   \text{ with }
   \P_{\bmf(\nnu)-\1} := \bigotimes_{m=1}^M \P_{\mf(\nu_m)-1}
   \text{ and }
   \bmf(\nnu) := (\mf(\nu_1), \ldots, \mf(\nu_M)).
\]
We denote by $\Lagr_m^{\mf(\nu_m)} : C^0(\G_m;\X) \to \P_{\mf(\nu_m)-1}(\G_m;\X)$
the univariate Lagrange interpolation operator associated with the set of nodes $\Colpts_m^{\mf(\nu_m)} \subset \G_m$.
Setting $\Lagr_m^{0} = 0$ for all $m = 1,\ldots,M$, we define univariate detail operators
\[
   \Delta_m^{\mf(\nu_m)} : = \Lagr_m^{\mf(\nu_m)} - \Lagr_m^{\mf(\nu_m-1)}.
\]
Now, the sparse grid collocation operator associated with the sparse grid $\Colpts_{\indset_\bullet}$ is defined~as
\begin{equation} \label{eq:S:def}
   S_{\bullet} = S_{\indset_\bullet} := \sum\limits_{\nnu \in \indset_\bullet} \Delta^{\bmf(\nnu)},
\end{equation}
where $\Delta^{\bmf(\nnu)} := \bigotimes_{m=1}^{M} \Delta_m^{\mf(\nu_m)}$ denotes the hierarchical surplus operator.

The operator $S_{\indset_\bullet}$ can be written also as a linear combination of tensor products of univariate
Lagrange interpolation operators as follows:
\begin{equation} \label{eq:S:implement}
   S_{\indset_\bullet} = \sum\limits_{\nnu \in \indset_\bullet} c_{\nnu} \bigotimes_{m=1}^{M} \Lagr_m^{\mf(\nu_m)}
   \text{\ \ with\ \ }
   c_{\nnu} := \sum_{\substack{\jj \in \{0, 1\}^M \\ (\nnu + \jj) \in \indset_\bullet}} (-1)^{|\jj|_1}.
\end{equation}
\rbl{This representation generates an efficient implementation of $S_{\indset_\bullet}$.} 
Furthermore, the nestedness of univariate node sets and the monotonicity of the index set $\indset_\bullet$ 
imply the interpolation property for the operator $S_{\indset_\bullet}$, i.e.,
\begin{equation} \label{eq:interpolation}
   S_{\indset_\bullet}: C^0(\G;\X) \to \P_{\indset_\bullet}(\G;\X)
   \text{\; is such that \;}
   S_{\indset_\bullet} v(\z) = v(\z)\ \; \forall\,\z \in \Colpts_{\indset_\bullet}.
\end{equation}
Therefore, the SC-FEM solution defined by~\eqref{eq:scfem:sol} can be written as
\begin{equation} \label{eq:scfem:sol:S-form}
   \scsol(x,\y) = S_{\bullet} U_\bullet (x,\y) = \sum\limits_{\z \in \Colpts_\bullet} u_{\bullet \z}(x) \LagrBasis{\bullet \z}{}(\y)
\end{equation}
with a function $U_{\bullet} : \G \to \SS^1_0 \Big( \bigoplus_{\z \in \Colpts_\bullet} \TT_{\bullet \z} \Big)$
satisfying $U_{\bullet}(\z) = u_{\bullet \z}$ for all $\z \in \Colpts_\bullet$;
here, $\bigoplus_{\z \in \Colpts_\bullet} \TT_{\bullet \z}$ denotes the \rbl{\it overlay} of the meshes
$\TT_{\bullet \z}$, $\z \in \Colpts_\bullet$ (\rbl{in other words}, their coarsest common refinement).

The enhancement of the parametric component of the SC-FEM approximation given by~\eqref{eq:scfem:sol:S-form}
is done \rbl{by} enriching the index set $\indset_\bullet$ (and, hence, expanding the
set $\Colpts_\bullet$ of collocation points).
To that end, for a given index set $\indset_\bullet$, we introduce the margin 
\begin{equation} \label{eq:margin}
   \marg_{\bullet} = \marg({\indset_\bullet}) :=
   \{ \nnu \in \N^{M} \backslash \indset_\bullet :
   \nnu - \eeps_m \in \indset_\bullet \text{ for some } m \in \{1, \dots, M\}\}
\end{equation}
and the reduced margin
\begin{equation} \label{eq:reduced:margin}
   \rmarg_{\bullet} = \rmarg({\indset_\bullet}) :=
   \{ \nnu \in \marg(\indset_\bullet) :
   \nnu - \eeps_m \in \indset_\bullet \text{ for all } m = 1, \dots, M \text{ such that } \nu_m > 1 \}.
\end{equation}
Note that for a monotone $\indset_\bullet$ and for any subset  \rbl{of marked indices} 
$\text{M}_\bullet \subseteq \rmarg_\bullet$, the index set $\indset_\bullet \cup \text{M}_\bullet$ is also~monotone.

\section{Hierarchical a posteriori error estimation and adaptivity} \label{sec:error:estimate}

\rbl{In the sequel, we define $\norm{\cdot}{\X} := \norm{\nabla\cdot}{L^2(D)}$ and let
$\norm{\cdot}{}$ denote the norm in $\V = L^p_{\pi}(\G,\X)$ for a fixed $1 \le p \le \infty$.}
We \rbl{will use a hierarchical construction}
(see, e.g.,~\cite[Chapter~5]{AinsworthOden00}) to derive a reliable a posteriori estimate
for the discretization error $u - \scsol = u - S_{\bullet} U_\bullet \in \V$.
To that end, we denote by $\scsolhat$
an enhanced SC-FEM approximation that reduces the discretization error, i.e.,
\begin{equation} \label{eq:saturation}
   \norm{u - \scsolhat}{} \le q_{\rm sat} \norm{u - \scsol}{} 
\end{equation}
with some constant $q_{\rm sat} \in (0,1)$ \rev{that is independent of discretization parameters}.
Then, by using the triangle inequality, we~obtain
\begin{equation} \label{eq:estimate:1}
   \norm{u - \scsol}{} \le \big(1 - q_{\rm sat}\big)^{-1} \, \norm{\scsolhat - \scsol}{}.
\end{equation}
\abrevx{We consider the following enhanced solution}
\begin{equation} \label{eq:scfem:sol:enhanced}
\abrevx{
   \scsolhat := S_\bullet \widehat U_\bullet + \left( \widehat S_\bullet \widetilde U_\bullet - S_\bullet U_\bullet\right),}
\end{equation}
\rbl{where}
\begin{equation} \label{eq:hat:U}
   \widehat U_{\bullet} : \G \to \SS^1_0 \Big( \bigoplus_{\z \in \Colpts_\bullet} \widehat\TT_{\bullet \z} \Big)
   \text{ is such that }
   \widehat U_{\bullet}(\z) = \widehat u_{\bullet \z} 
   \ \ \forall\, \z \in \Colpts_\bullet,
\end{equation}
\begin{equation} \label{eq:hat:S}
   \abrevx{
   \widehat S_{\bullet} = S_{\widehat\indset_\bullet} := \sum\limits_{\nnu \in \widehat\indset_\bullet} \Delta^{\bmf(\nnu)},}
\end{equation}
\abrevx{and}
\begin{equation} \label{eq:tilde:U}
   \widetilde U_{\bullet} :
   \G \to \SS^1_0 
 \rbl{   \Big( \Big[ \bigoplus_{\z \in  \widehat\Colpts_\bullet} \TT_{\bullet \z} \Big]  \Big)}
   \text{ \abrevx{is} \rbl{defined so that} }
   \widetilde U_{\bullet}(\z) = \widetilde u_{\bullet \z} = 
   \begin{cases}
      u_{\bullet \z} 
      & \!\!\forall\, \z \,{\in}\, \Colpts_\bullet,\\
      u_{0 \z} 
      & \!\!\forall\, \z \,{\in}\, \widehat\Colpts_\bullet \setminus \Colpts_\bullet
   \end{cases} \abrevx{.}
\end{equation}
\rbl{Here, to retain generality,} $\widehat\indset_\bullet$ is \rbl{any} monotone index set that contains $\indset_\bullet$
(e.g., $\widehat\indset_\bullet = \indset_\bullet \cup \rmarg_\bullet$)\abrevx{, and}
$\widehat\Colpts_\bullet$ is the set of collocation points generated by the index set $\widehat\indset_\bullet$.

\rbl{A subtle feature of the construction \eqref{eq:tilde:U}
 is the identification of
  $u_{0 \z} \in \SS^1_0(\rbl{\TT_{0\z}})$ as} the Galerkin approximation on 
 \rbl{a suitable} (coarse) mesh~$\rbl{\TT_{0\z}}$ using the coefficient $a$ sampled at
 \rbl{a 
 \abrevx{new} collocation point} $\z \in \widehat\Colpts_\bullet \setminus \Colpts_\bullet$.
\rbl{The construction of 
  sample-specific  meshes for new collocation points will be discussed in detail in part~II of this work.}

\rbl{To summarise,} \abrevx{the definition of $\scsolhat$ in~\eqref{eq:scfem:sol:enhanced} is based on}
\rbl{two} enhanced (multilevel) SC-FEM approximations; \rbl{namely,}
\begin{itemize}
\item[(i)]
$S_\bullet \widehat U_\bullet$ \rbl{is determined by} the same set $\Colpts_\bullet$ of collocation points 
as the SC-FEM solution $\scsol$
but employs the enhanced Galerkin approximations $\widehat u_{\bullet \z} \in \widehat \X_{\bullet \z}$,
\rbl{and}
\item[(ii)]
$\widehat S_\bullet \widetilde U_\bullet$  \rbl{is determined by} the same
 Galerkin approximations $u_{\bullet \z} \in \X_{\bullet \z}$ as $\scsol$
at each collocation point $\z \in \Colpts_\bullet$ \rbl{in combination with
\abrevx{(coarse) mesh Galerkin approximations} $u_{0 \z} \in \SS_0^1(\rev{\TT_{0\z}})$
at all new collocation points} $\z \in \widehat\Colpts_\bullet \setminus \Colpts_\bullet$.
\abrevx{The term
$(\widehat S_\bullet \widetilde U_\bullet - S_\bullet U_\bullet)$ in~\eqref{eq:scfem:sol:enhanced}
is the hierarchical surplus associated with the enhanced approximation $\widehat S_\bullet \widetilde U_\bullet$.}
\end{itemize}

\rev{
\begin{remark} \label{rem:scfem:sol:enhanced:alt}
An alternative construction of a function $\scsolhat$  is given by
\begin{equation} \label{eq:scfem:sol:enhanced:alt}
   \scsolhat := S_\bullet \widehat U_\bullet + 
   \left( \widehat S_\bullet \widetilde U_{\bullet,0} - S_\bullet U_{\bullet,0}\right),
\end{equation}
where $\widehat U_{\bullet}$ and $\widehat S_{\bullet}$ are defined as before by~\eqref{eq:hat:U} 
and~\eqref{eq:hat:S}, respectively,
\begin{equation*} 
   \widetilde U_{\bullet,0} :
   \G \to \SS^1_0(\TT_0) 
   \text{ is such that }
   \widetilde U_{\bullet,0}(\z') = u_{0 \z'}\ \forall\, \z' \in \widehat\Colpts_\bullet,
\end{equation*}
and
\begin{equation*}
   U_{\bullet,0} :
   \G \to \SS^1_0(\TT_0) 
   \text{ is such that }
   U_{\bullet,0}(\z) = u_{0 \z}\ \forall\, \z \in \Colpts_\bullet.
\end{equation*}
The advantage of this construction compared to that  in~\eqref{eq:scfem:sol:enhanced}
 is in the ease of implementation of the parametric enhancement
$\widehat S_\bullet \widetilde U_{\bullet,0} - S_\bullet U_{\bullet,0}$, as the involved Galerkin approximations
$u_{0 \z'}$ ($\z' \in \widehat\Colpts_\bullet$) and $u_{0 \z}$ ($\z \in \Colpts_\bullet$) are all computed
on the coarsest finite element mesh $\TT_0$.
\end{remark}
}

We assume that $\scsolhat$ defined by \eqref{eq:scfem:sol:enhanced} satisfies the saturation 
property~\eqref{eq:saturation}.  Therefore, by using~\eqref{eq:scfem:sol:S-form},~\eqref{eq:scfem:sol:enhanced} 
and the triangle inequality, we derive from~\eqref{eq:estimate:1}
\begin{equation} \label{eq:estimate:2}
   \norm{u - \scsol}{} \le \frac{1}{\abrevx{1 - q_{\rm sat}}} \,
                                      \Big(
                                      \norm{S_\bullet(\widehat U_\bullet - U_\bullet)}{} +
                                      \norm{\widehat S_\bullet \widetilde U_\bullet - S_\bullet U_\bullet}{}
                                      \Big).
\end{equation}
The two norms on the right hand side of~\eqref{eq:estimate:2} can be seen as the spatial and the parametric components
of an a posteriori estimate for the discretization error.
We will denote these spatial and parametric error estimates as
\begin{equation} \label{eq:two:estimates}
   \mu_\bullet := \norm{S_\bullet(\widehat U_\bullet - U_\bullet)}{}\quad
   \text{and}\quad
   \tau_\bullet := \norm{\widehat S_\bullet \widetilde U_\bullet - S_\bullet U_\bullet}{},
\end{equation}
respectively.
\rbl{The components of the  error estimate are discussed in more detail  below}.

\subsection{Spatial error estimate and spatial error indicators} \label{sec:error:spatial}

For the spatial error estimate $\mu_\bullet$, 
\rbl{we can estimate component error contributions}
 using the triangle inequality
\begin{align}
   \mu_\bullet = \norm{S_\bullet (\widehat U_\bullet - U_\bullet)}{}
   & \reff{eq:scfem:sol:S-form} = \norm[\bigg]{\sum\limits_{\z \in \Colpts_\bullet}
    (\widehat u_{\bullet \z} - u_{\bullet \z})\, \LagrBasis{\bullet \z}{}}{}
   \le \sum\limits_{\z \in \Colpts_\bullet} \norm{\widehat u_{\bullet \z} - 
   u_{\bullet \z}}{\X} \, \norm{\LagrBasis{\bullet \z}{}}{L^p_\pi(\G)}.
   \label{eq:estimate:3}
\end{align}
\rbl{The crude  bound in~\eqref{eq:estimate:3} pinpoints the inbuilt advantage of  stochastic Galerkin 
approximation over SC-FEM approximation in the context of solving elliptic PDEs with random 
data. Numerical experiments confirm that the componentwise  bound gives nonrobust over-estimation 
of the error as the number of parameters is  increased.}

Despite the inaccuracy, we will demonstrate later that the componentwise bound can still be 
employed to define an error {\it indicator} that can be used to drive a reliable adaptive refinement process.
While the norm $\norm{\widehat u_{\bullet \z} - u_{\bullet \z}}{\X}$ is computable, its evaluation
for each $\z \in \Colpts_\bullet$ requires computation of the enhanced Galerkin approximation 
$\widehat u_{\bullet \z}$.
It is computationally more efficient to \emph{estimate}  these error components using  hierarchical error indicators.
We review two \rev{possible approaches} below.

\smallskip
{\bf Spatial hierarchical error \rbl{indicator} I.}
For each $\z \in \Colpts_\bullet$, we define $\mu_{\bullet \z} := \norm{e_{\bullet \z}}{\X}$,
where $e_{\bullet \z} \in \Y_{\bullet \z}$ satisfies
(recall that $\widehat \X_{\bullet \z} = \X_{\bullet \z} \oplus \Y_{\bullet \z}$)
\begin{equation} \label{eq:hierar:estimator}
   \int_D \nabla e_{\bullet \z}(x) \cdot \nabla v(x) \dx =
   \int_D f(x) v(x) \dx -
   \int_D a(x, \z) \nabla u_{\bullet \z}(x) \cdot \nabla v(x) \dx \quad \forall v \in \Y_{\bullet \z}.
\end{equation}
Then, \rbl{following  the construction in~\cite{banksmith93} and using the norm equivalence 
in~\eqref{eq:norm:equiv} the following estimate holds}:
\begin{equation} \label{eq:estimate:4}
   \norm{\widehat u_{\bullet \z} - u_{\bullet \z}}{\X} \le
   a_{\min}^{-1}\, (1 - \g_\z^2)^{-1/2}\, \mu_{\bullet \z}, 
   \end{equation}
where $\g_\z \in [0,1)$ is the constant in the strengthened Cauchy--Schwarz inequality~\eqref{eq:sCS}
for the subspaces $\X_{\bullet \z}$, $\Y_{\bullet \z}$.
\rbl{Moreover, from~\eqref{eq:estimate:3} we get the crude overestimate}
\begin{equation} \label{eq:estimate:5}
   \mu_\bullet = \norm{S_\bullet (\widehat U_\bullet - U_\bullet)}{} \le
   a_{\min}^{-1}\, (1 - \g^2)^{-1/2}\,
   \sum\limits_{\z \in \Colpts_\bullet} \mu_{\bullet \z} \, \norm{\LagrBasis{\bullet \z}{}}{L^p_\pi(\G)},
\end{equation}
where $\g := \max\{\g_\z : \z \in \Colpts_\bullet\}$.
Note that the constant $a_{\min}^{-1}$ in~\eqref{eq:estimate:5} is the same as in~\cite{GuignardN18}.
\begin{remark}
\rbl{It is important that $\g_{\z}$ in~\eqref{eq:estimate:4} is independent 
of the coefficient sample $a(\cdot,\z)$. Noting that $\g$ only 
depends on the subspaces $\X_{\bullet \z}$, $\Y_{\bullet \z}$ and that 
 all underlying triangulations are generated from the same coarse mesh $\TT_0$ by applying NVB refinement(s),
 it \abrevx{is} eminently plausible that there exists a uniform upper bound for $\g$ for all nested sets of
  collocation points  generated when running an adaptive algorithm.
An alternative construction would include the coefficient 
sample $a(\cdot,\z)$ in definition of the problem \eqref{eq:hierar:estimator}. If this was done then
the constants $\tilde\g_\z$ would depend on the samples of the coefficient $a(\cdot,\z)$ at 
collocation points in which case it is not so obvious that there exists a uniform upper bound for $\g$.}
\end{remark}
While local error indicators are not explicitly defined in the above construction, they
can be easily derived from the computed error estimator~$e_{\bullet \z}$.
For example, for each element $T \in \TT_{\bullet \z}$, the (spatial) error indicator associated 
with $T$ is given by
\[
   \mu_{\bullet \z}(T) := \norm{ e_{\bullet \z} \vert_{T} }{\X}.
\]
Alternatively, the (spatial) error indicators $\mu_{\bullet \z}(\xi)$ associated  with interior 
edge midpoints $\xi \in \NN_{\bullet \z}^+$ are given by components of the solution vector to the linear system 
stemming from the discrete formulation~\eqref{eq:hierar:estimator}.

\smallskip
{\bf Spatial hierarchical error \rbl{indicator} II.} 
Recall that 
$\Y_{\bullet \z} := \hull\{ \widehat\varphi_{\bullet \z,\xi} : \xi \in \NN_{\bullet \z}^+ \}$ (see~\eqref{eq:Y:space}).
For each $\z \in \Colpts_\bullet$, we \rbl{can} define the two-level error indicators
associated  with interior edge midpoints:
\begin{equation} \label{eq:2level:indicator}
   \mu_{\bullet \z}(\xi) :=
   \frac{\big| (f,\widehat\varphi_{\bullet \z,\xi})_{L^2(D)} - (a(\cdot,\z) \nabla u_{\bullet \z}, \nabla \widehat\varphi_{\bullet \z,\xi})_{L^2(D)} \big|}
           {\norm{\widehat\varphi_{\bullet \z,\xi}}{\X}}\quad
   \text{for all } \xi \in \NN_{\bullet \z}^+.
\end{equation}
These indicators \rbl{can then be} combined to produce the two-level \rbl{indicator}
\begin{equation}
   \mu^2_{\bullet \z} := \sum\limits_{\xi \in \NN_{\bullet \z}^+} \mu^2_{\bullet \z}(\xi)
\end{equation}
that satisfies (see~\cite{msw98, ms99})
\begin{equation} \label{eq:estimate:6}
   \norm[\big]{a^{1/2}(\cdot,\z) \, \nabla(\widehat u_{\bullet \z} - u_{\bullet \z})}{L^2(D)} \le
   a_{\min}^{-1/2}\, C_{\rm est}\, \mu_{\bullet \z},
\end{equation}
where $a_{\min}$ is the constant in~\eqref{eq:norm:equiv} and \rbl{$C_{\rm est}$} is a generic constant
that only depends on the \abrevx{coarse mesh $\TT_0$}.

Using again the norm equivalence in~\eqref{eq:norm:equiv}, \rbl{we get the following 
estimate of the error\abrevx{:}}
\begin{align}
   \mu_\bullet = \norm{S_\bullet (\widehat U_\bullet - U_\bullet)}{} & \le
   a_{\min}^{-1/2}\,
   \sum\limits_{\z \in \Colpts_\bullet}
   \norm[\big]{a^{1/2}(\cdot,\z) \, \nabla(\widehat u_{\bullet \z} - u_{\bullet \z})}{L^2(D)} \, \norm{\LagrBasis{\bullet \z}{}}{L^p_\pi(\G)}
   \nonumber\\[4pt]
   & \le
   a_{\min}^{-1}\, C_{\rm est}\,
   \sum\limits_{\z \in \Colpts_\bullet}
   \mu_{\bullet \z}\, \norm{\LagrBasis{\bullet \z}{}}{L^p_\pi(\G)}.
   \label{eq:estimate:7}
\end{align}

\rbl{The advantage of hierarchical error estimators over the residual estimators
 discussed by Guignard \& Nobile in~\cite{GuignardN18} is that they provide
  information about potential error reduction associated with local refinement in 
  space or \rev{with} enhancement of the parametric approximation, see~\cite{bps14,bs16}.
  They also provide a more natural starting point for a rigorous convergence analysis
  of adaptive strategies, see~\cite{bprr18}.}

\subsection{Parametric error estimate and parametric error indicators} \label{sec:error:parametric}

\rbl{We now focus} on the parametric error estimate $\tau_\bullet$ in~\eqref{eq:two:estimates}.
\abrev{First, recalling the definitions of the operators $S_\bullet$ and $\widehat S_\bullet$
in~\eqref{eq:S:def} and~\eqref{eq:hat:S}, respectively, we find that
\begin{align} \label{eq:estimate:7:1}
   \widehat S_\bullet \widetilde U_\bullet - S_\bullet U_\bullet & =
   \sum\limits_{\nnu \in \widehat\indset_\bullet} \Delta^{\bmf(\nnu)} \widetilde U_\bullet - 
   \sum\limits_{\nnu \in \indset_\bullet} \Delta^{\bmf(\nnu)} U_\bullet =
   \nonumber
   \\[4pt]
   & =
   \sum\limits_{\nnu \in \indset_\bullet} \Delta^{\bmf(\nnu)} (\widetilde U_\bullet - U_\bullet) +
   \sum\limits_{\nnu \in \widehat\indset_\bullet \setminus \indset_\bullet} \Delta^{\bmf(\nnu)} \widetilde U_\bullet
   \reff{eq:tilde:U}=
   \sum\limits_{\nnu \in \widehat\indset_\bullet \setminus \indset_\bullet} \Delta^{\bmf(\nnu)} \widetilde U_\bullet.
\end{align}
On the other hand,}
thanks to $\indset_\bullet$ and $\widehat\indset_\bullet$ being monotone, we \abrev{can} write 
\begin{align*}
   \big( \widehat S_\bullet \widetilde U_\bullet - S_\bullet U_\bullet \big)(x,\y) & \refp{eq:tilde:U}=
   \sum\limits_{\z \in \widehat\Colpts_\bullet} \widetilde u_{\bullet \z}(x) \LagrBasisHat{\bullet \z}{}(\y) -
   \sum\limits_{\z \in \Colpts_\bullet} u_{\bullet \z}(x) \LagrBasis{\bullet \z}{}(\y)
   \\[4pt]
   & \reff{eq:tilde:U}=
   \sum\limits_{\z \in \Colpts_\bullet} u_{\bullet \z}(x) \big( \LagrBasisHat{\bullet \z}{}(\y) - \LagrBasis{\bullet \z}{}(\y) \big) +
   \sum\limits_{\z \in \widehat\Colpts_\bullet \setminus \Colpts_\bullet} u_{0 \z}(x) \LagrBasisHat{\bullet \z}{}(\y),
\end{align*}
where $\LagrBasisHat{\bullet \z}{}(\y) = \LagrBasis{\z}{\widehat\Colpts_\bullet}(\y)$ denotes
the Lagrange polynomial basis function associated with the point $\z \in \widehat\Colpts_\bullet$ and satisfying
$\LagrBasisHat{\bullet \z}{}(\z') = \delta_{\z\z'}$ for any $\z,\,\z' \in \widehat\Colpts_\bullet$.
Note that for any $\z' \in \widehat\Colpts_\bullet$ there holds
\begin{align*}
   \big( \widehat S_\bullet \widetilde U_\bullet - S_\bullet U_\bullet \big)(x,\z') & =
   \begin{cases}
      \sum\limits_{\z \in \Colpts_\bullet} u_{\bullet \z}(x) \big( \delta_{\z \z'} - \delta_{\z \z'} \big) +
      \sum\limits_{\z \in \widehat\Colpts_\bullet \setminus \Colpts_\bullet} u_{0 \z}(x) \delta_{\z \z'}
      & \text{if } \z' \in \Colpts_\bullet,
      \\[18pt]
      \sum\limits_{\z \in \Colpts_\bullet} u_{\bullet \z}(x) \big( \delta_{\z \z'} - \LagrBasis{\bullet \z}{}(\z') \big) +
      \sum\limits_{\z \in \widehat\Colpts_\bullet \setminus \Colpts_\bullet} u_{0 \z}(x) \delta_{\z \z'}
      & \text{if } \z' \in \widehat\Colpts_\bullet \setminus \Colpts_\bullet
   \end{cases}
   \\[4pt]
   & =
   \begin{cases}
      0 & \text{if } \z' \in \Colpts_\bullet,
      \\[4pt]
      -\scsol(x,\z') + u_{0 \z'}(x) & \text{if } \z' \in \widehat\Colpts_\bullet \setminus \Colpts_\bullet.
   \end{cases}
\end{align*}
Therefore,
\begin{equation} \label{eq:estimate:8}
   \tau_\bullet = \norm[\big]{\widehat S_\bullet \widetilde U_\bullet - S_\bullet U_\bullet}{} =
   \norm[\bigg]{\sum\limits_{\z \in \widehat\Colpts_\bullet \setminus \Colpts_\bullet}
                        \big( u_{0 \z} - \scsol(\cdot,\z) \big) \LagrBasisHat{\bullet \z}{}}{}.
\end{equation}
The parametric estimate $\tau_\bullet$ is thus computable;
 \rbl{calculating it} requires extra PDE solves on \rbl{coarse meshes~$\TT_{0\z}$ for a small 
number} of collocation points $\z \in \widehat\Colpts_\bullet \setminus \Colpts_\bullet$.

\abrev{The natural parametric \abrevx{error indicators}
\rbl{associated with} \eqref{eq:estimate:7:1} are given by
\begin{equation} \label{eq:param:indicators}
   \tau_{\bullet \nnu} :=
   \norm[\big]{\Delta^{\bmf(\nnu)} \widetilde U_\bullet}{},\qquad
   \nnu \in \widehat\indset_\bullet \setminus \indset_\bullet.
\end{equation}
}
\begin{remark} \label{rem:indicator:parametric}
If the enriched index set $\widehat\indset_\bullet$ is obtained using the reduced margin of $\indset_\bullet$,
i.e., $\widehat\indset_\bullet \setminus \indset_\bullet = \rev{\rmarg}(\indset_\bullet)$,
then the collocation points in the set $\widehat\Colpts_\bullet \setminus \Colpts_\bullet$ can be grouped 
together according to the `generating' multi-index $\nnu \in \widehat\indset_\bullet \setminus \indset_\bullet$ 
such that
\[
   \widehat\Colpts_\bullet \setminus \Colpts_\bullet =
   \mathop{\cup}\limits_{\nnu \in \widehat\indset_\bullet \setminus \indset_\bullet} \widetilde\Colpts_{\bullet \nnu}
   \quad
   \text{and}
   \quad
   \widetilde\Colpts_{\bullet \nnu} \cap \widetilde\Colpts_{\bullet \nnu'} = \emptyset\ \ 
   \forall\,\nnu, \nnu' \in \widehat\indset_\bullet \setminus \indset_\bullet,\; \nnu \not= \nnu'.
\]
In this case, we conclude from~\eqref{eq:estimate:7:1} and~\eqref{eq:estimate:8} that
\[
   \tau_\bullet = \norm[\big]{\widehat S_\bullet \widetilde U_\bullet - S_\bullet U_\bullet}{} =
   \norm[\bigg]{\sum\limits_{\nnu \in \widehat\indset_\bullet \setminus \indset_\bullet}\;
                        \sum\limits_{\z \in \widetilde\Colpts_{\bullet \nnu}}
                        \big( u_{0 \z} - \scsol(\cdot,\z) \big) \LagrBasisHat{\bullet \z}{}}{}\abrevx{,}
\]
and the associated parametric \abrevx{error indicators} are given by
\begin{equation} \label{eq:param:indicators:1}
   \widetilde\tau_{\bullet \nnu} =
   \sum\limits_{\z \in \widetilde\Colpts_{\bullet \nnu}}
   \norm{u_{0 \z} - \scsol(\cdot,\z)}{\X} \, \norm{\LagrBasisHat{\bullet \z}{}}{L_\pi^p(\G)},\qquad
   \nnu \in \widehat\indset_\bullet \setminus \indset_\bullet.
\end{equation}
Note that for linearly growing sets of Leja points, one has $\tau_{\bullet \nnu} = \widetilde\tau_{\bullet \nnu}$,
since for each $\nnu \in \widehat\indset_\bullet \setminus \indset_\bullet$, the set $\widetilde\Colpts_{\bullet \nnu}$ consists of
a single point $\z_{\nnu} \in \widehat\Colpts_\bullet \setminus \Colpts_\bullet$.
\end{remark}

\rev{
\begin{remark} \label{rem:indicator:parametric:alt}
All the arguments in this section extend trivially to the parametric error estimate and parametric error indicators
derived from the alternative construction~\eqref{eq:scfem:sol:enhanced:alt} of the enhanced 
SC-FEM solution~$\scsolhat$.
\end{remark}
}

\subsection{Adaptivity} \label{sec:adaptivity}
\rbl{A general multilevel SC-FEM adaptive algorithm is presented below. 
There are two features worth noting at the outset.
 First,  the refinement of finite element approximations and the enrichment of the set of 
 collocation points are driven by \abrevx{the error indicators} $ \mu_{\bullet \z}$  and $ \tau_{\bullet \abrevx{\nnu}}$
 discussed above. 
 Second,   the  error estimates  $\mu_\bullet$  and  $\tau_\bullet$ in \eqref{eq:two:estimates} only need to be 
 calculated periodically;  their combination is required for reliable termination of the adaptive process and  to
 provide reassurance that the SC-FEM error is decreasing at an acceptable rate.}
 
\begin{algorithm} \label{algorithm}
{\bfseries Input:}
$\indset_0 = \{ \1 \}$; 
\rbl{$\TT_{0 \z}$} for all $\z \in \widehat\Colpts_0 := \Colpts_{\indset_0 \cup \rmarg(\indset_0)}$;
marking criterion.\\
Set the iteration counter $\ell := 0$,  the output counter $k$ and the \rev{error} tolerance.
\begin{itemize}
\item[\rm(i)] 
Compute Galerkin approximations $\big\{ u_{\ell \z} \in \X_{\ell \z} : \z \in \abrevx{\widehat\Colpts_\ell} \big\}$ 
by solving~\eqref{eq:sample:fem}.
\item[\rm(ii)] 
Compute the spatial \abrevx{error indicators} $\big\{ \mu_{\ell\z}: \z \in \Colpts_\rev{\ell} \big\}$
\rev{by solving~\eqref{eq:hierar:estimator}}
(or the indicators $\big\{ \mu_{\ell \z}(\xi) : \z \in \Colpts_\ell,\ \xi \in \NN_{\ell \rev{\z}}^{+} \big\}$ 
given by~\eqref{eq:2level:indicator}).
\item[\rm(iii)] 
Compute the parametric \abrevx{error indicators}
\rbl{$\big\{ \tau_{\ell \nnu} : \nnu \in \widehat\indset_\ell \setminus \indset_\ell \big\}$
given by~\eqref{eq:param:indicators}
(or the indicators $\big\{ \widetilde\tau_{\ell \nnu} : \nnu \in \widehat\indset_\ell \setminus \indset_\ell \big\}$
given by~\eqref{eq:param:indicators:1}).}
\item[\rm(iv)] 
Use a marking criterion to determine $\MM_{\ell \z} \subseteq \NN_{\ell \z}^+$ for all $\z \in \Colpts_\ell$ and
\abrev{$\markindset_\ell \subseteq \widehat\indset_\ell \setminus \indset_\ell$}.
\item[\rm(v)] For all $\z \in \Colpts_\ell$, set $\TT_{(\ell+1) \z} := \refine(\TT_{\ell \z},\MM_{\ell \z})$.
\item[\rm(vi)] Set \abrev{$\indset_{\ell+1} := \indset_\ell \cup \markindset_\ell$,}
$\widehat\Colpts_{\ell+1} := \Colpts_{\indset_{\ell+1} \cup \rmarg(\indset_{\ell+1})}$, and
\rbl{construct  $\TT_{(\ell+1) \z} := \rev{\TT_{0 \z}}$} for all $\z \in \widehat\Colpts_{\ell+1} \setminus \Colpts_\rev{\ell}$.
\item[\rm(vii)] 
If $\ell = j k , j\in \N$, compute the spatial and parametric  {error estimates}  $\mu_\ell$  and  $\tau_\ell$ 
and exit if $\mu_\ell + \tau_\ell < {\tt error tolerance}$.
\item[\rm(viii)] Increase the counter $\ell \mapsto \ell+1$ and goto {\rm(i)}.
\end{itemize}
{\bfseries Output:} For some specific  $\rbl{\ell_*=jk} \in \N$,
the algorithm returns the multilevel SC-FEM approximation~$u_{\rbl{\ell_*}}^{\rm SC}$
computed via~\eqref{eq:scfem:sol} from Galerkin approximations 
$\big\{ u_{\rbl{{\ell_*}} \z} \in \X_{\rbl{{\ell_*}} \z} : \z \in \Colpts_\ell \big\}$
together with a corresponding \rev{error estimate} \abrevx{$\mu_{\ell_*} + \tau_{\ell_*}$}.
\end{algorithm}

A \rev{simple} marking strategy for step~(iv) of Algorithm~\ref{algorithm} is specified next.

\medskip
\noindent {\bf Marking criterion.}
\textbf{Input:}
error \rbl{reduction}  indicators $\big\{ \mu_{\ell \z}(\xi) : \z \in \Colpts_\ell,\ \xi \in \NN_{\ell \z}^+ \big\}$,
$\big\{ \mu_{\ell \z} = \big( \sum_{\xi \in \NN_{\ell \z}^+} \mu^2_{\ell \z}(\xi) \big)^{1/2}: \z \in \Colpts_\ell \big\}$,
and
$\big\{ \tau_{\ell \nnu} : \nnu \in \widehat\indset_\ell \setminus \indset_\ell \big\}$;
marking parameters $0 < \theta_\X, \theta_\Colpts \le 1$ and $\vartheta > 0$.
\begin{itemize}
\item[$\bullet$]
If \
$\sum_{\z \in \Colpts_\ell} \mu_{\ell \z} \norm{L_{\ell \z}}{L^p_{\pi}(\G)}
  \ge \vartheta \sum_{\abrev{\nnu \in \widehat\indset_\ell \setminus \indset_\ell}} \abrev{\tau_{\ell \nnu}}$,
then proceed as follows:
\begin{itemize}
\item[$\circ$]
set $\abrev{\markindset_\ell} := \emptyset$
\item[$\circ$]
for each $\z \in \Colpts_\ell$,
determine $\MM_{\ell \z} \subseteq \NN_{\ell \z}^+$ of minimal cardinality such that
\begin{equation} \label{eq:doerfler:separate1}
 \theta_\X \, \mu_{\ell \z}^2 
 \le \sum_{\xi \in \MM_{\ell \z}} \mu^2_{\ell \z}(\xi).
\end{equation}
\end{itemize}
\item[$\bullet$]
Otherwise, if \ 
$\sum_{\z \in \Colpts_\ell} \mu_{\ell \z} \norm{L_{\ell \z}}{L^p_{\pi}(\G)}
  < \vartheta \sum_{\abrev{\nnu \in \widehat\indset_\ell \setminus \indset_\ell}} \abrev{\tau_{\ell \nnu}}$,
then proceed as follows:
\begin{itemize}
\item[$\circ$]
set $\MM_{\ell \z} := \emptyset$ for all $\z \in \Colpts_\ell$
\item[$\circ$]
determine \abrev{$\markindset_\ell \subseteq \widehat\indset_\ell \setminus \indset_\ell$}
of minimal cardinality such that
\begin{equation} \label{eq:doerfler:separate2}
 \theta_\Colpts \, \sum_{\abrev{\nnu \in \widehat\indset_\ell \setminus \indset_\ell}} \abrev{\tau_{\ell \nnu}} \le
 \sum_{\abrev{\nnu \in \markindset_\ell}} \abrev{\tau_{\ell \nnu}}.
\end{equation}
\end{itemize}
\end{itemize}
\textbf{Output:}
$\MM_{\ell \z} \subseteq \NN_{\ell \z}^+$ for all $\z \in \Colpts_\ell$ and
\abrev{$\markindset_\ell \subseteq \widehat\indset_\ell \setminus \indset_\ell$}.

\medskip
\rbl{The rationale for checking convergence periodically  in Algorithm~\ref{algorithm}
(rather than every iteration) is that  direct computation of the spatial and parametric 
error estimates in \eqref{eq:two:estimates} incurs a significant computational overhead. 
\abrevx{In particular}, the calculation of the spatial error estimate $\mu_\ell$ requires the solution of 
 the PDE on \rev{uniform refinements} of all \rev{meshes} associated with 
 collocation points  generated by the current index set. 
}

\rbl{
\abrevx{We set $p = 2$ when computing the norms in $\V = L^p_{\pi}(\G,\X)$ in practice.}
The only other detail needed to implement Algorithm~\ref{algorithm} 
\abrevx{is} the specification of  the starting \rev{meshes} $\TT_{0 \z}$  when introducing new collocation 
points\footnote{\rbl{Starting \rev{meshes} also need to be specified in the initialization phase.}} 
in step (vi).  This specification will be shown to  be crucially important  in part~II of
this work.}
In the standard {\it single-level} SC-FEM  setting discussed below,
\rev{the same current mesh $\TT_{\bullet}$ is assigned to
all new collocation points added in step~(vi) of Algorithm~\ref{algorithm}.
Accordingly,}
the meshes $\TT_{0\z}$ ($\z \in \widehat\Colpts_\bullet \setminus \Colpts_\bullet$)
are \rev{set to be} identical to the meshes
$\TT_{\bullet \z} = \TT_{\bullet}$ ($\z \in \Colpts_\bullet$).
Thus, the construction of an overlay mesh when computing $\widetilde u_{\bullet \z}$ in~\eqref{eq:tilde:U} is
not needed \rev{in the single-level setting}.
 
\section{Numerical experiments}\label{sec:results}

The \abrevx{numerical} results \abrevx{presented in this section} show that adaptive SC-FEM strategies 
are competitive  in terms of computational effort with single-level adaptive stochastic Galerkin \rev{(SG)}
approximation---certainly in the  context of the model problem that is 
the focus of this study. The results also provide  a basis for comparison with
\abrevx{multilevel} adaptive SC-FEM in part II.

The single-level refinement strategy that will be employed is the obvious and
natural   simplification of the multilevel strategy described in \S\ref{sec:adaptivity}. 
Thus, at each step $\ell$ of the process, 
we compute the \abrevx{error indicators} associated with the SC-FEM solution $u_{\ell\z}$
(steps (ii)--(iii) of Algorithm~\ref{algorithm}).
\abrevx{Specifically, we employ the spatial hierarchical error indicator~I
computed \rev{by solving~\eqref{eq:hierar:estimator}}
and the parametric error indicators given by~\eqref{eq:param:indicators}.}
\abrevx{The marking criterion listed in~\S\ref{sec:adaptivity} (we set $\vartheta =1$) then identifies the refinement type by comparing}
the (global) spatial \rev{error estimate}
$\bar\mu_{\ell} := \| \abrevx{\mu_{\ell \z} \norm{L_{\ell \z}}{L^p_{\pi}(\G)}} \|_{\ell_1}$ with
the  parametric \rev{error estimate} $\bar\tau_{\ell} := \| \tau_{\ell\abrevx{\nnu}} \|_{\ell_1}$.
\abrevx{Thus,} if $\bar\mu_\ell$
is less than $\bar\tau_\ell$ then we enforce  a parametric refinement: choosing a bigger index set but 
keeping the finite element space unchanged; otherwise, we effect a spatial refinement: 
choosing an enhanced finite element space but keeping the index set unchanged.
The marking strategy \abrevx{also} generates the refinement process 
with the qualification that to effect  a spatial refinement,
we use D{\" o}rfler marking \abrevx{with  marking threshold $\theta_\X$} to produce sets of marked elements 
from the (single) \rev{mesh} $\TT_\ell$. 
A refined triangulation $\TT_{\ell+1}$ can then be constructed
by refining the elements in the \abrevx{{\it union}} of these individual
\abrevx{sets $\MM_{\ell \z}$ ($\z \in \Colpts_\ell$) of marked elements}.

\subsection{Test case I:  affine coefficient data}\label{sec:affineresults}

We set $f = 1$ and look to solve the model problem on the square-shaped 
domain \rev{$D = (0, 1)^2$} with random field coefficient given by
\begin{align} \label{kl}
a(x, \y) = a_0(x) + \sum_{m = 1}^M a_m(x) \, y_m,\quad
   x \in D,\ \y \in \Gamma.
\end{align}
The  specific problem  we consider is taken  from~\cite{bs16}.
The parameters $y_m$ in \eqref{kl} are the images of uniformly 
distributed independent mean-zero random variables, so that $\pi_m = \pi_m(y_m)$ is
the associated probability measure on $\G_m = [-1,1]$.
The expansion coefficients $a_m$, $m\, \in \, \N_0$  are chosen
to represent planar Fourier modes of increasing total order.
Thus, we fix  $a_0(x) := 1$ and set
\begin{equation}
\label{diff_coeff_Fourier}
a_m(x) := \alpha_m \cos(2\pi\beta_1(m)\,x_1) \cos(2\pi\beta_2(m)\,x_2),\  x=(x_1,x_2) 
\in (0,1) \times (0,1).
\end{equation}
The modes are ordered so that for any $m \in \N$,
\begin{equation}
  \beta_1(m) = m - k(m)(k(m)+1)/2\ \ \hbox{and}\ \ \beta_2(m) =k(m) - \beta_1(m)
\end{equation}
with $k(m) = \lfloor -1/2 + \sqrt{1/4+2m}\rfloor$ and the amplitude coefficients are 
constructed so that
$\alpha_m = \bar\alpha m^{-2}$  with $ \bar\alpha = 0.547$.
(This is referred to as the {\it slow decay case} in~\cite{bs16}.)
The  precise definition  of the amplitude coefficients ensures that 
the requirement \eqref{eq:amin:amax} is valid so the test problem 
is well posed in the sense discussed in section~\ref{sec:problem}.

\begin{figure}[!th]
\centering
\includegraphics[width = 0.75\textwidth]{{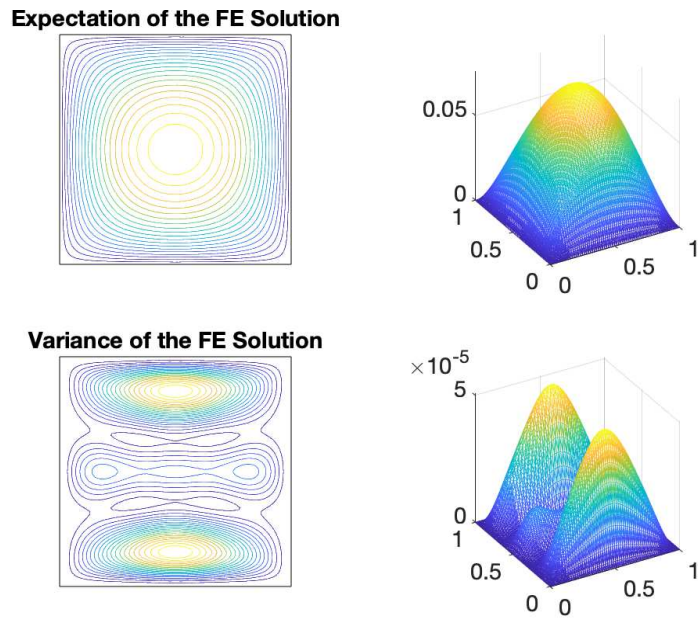}}
\caption{Reference solution for test case I.}
\label{fig:sc5.1sol}
\end{figure}

A reference solution to this problem is illustrated in Fig.~\ref{fig:sc5.1sol}. This solution 
was generated by running the adaptive algorithm with  $M$ set to \abrevx{4},  the marking
parameters \rev{$\theta_\X$} and \rev{$\theta_\Colpts$} both set to {\tt 0.3} and with 
$\tt{error tolerance}$ set to \abrevx{\tt 6e-3}.  This  tolerance was satisfied after 
 \abrevx{20} spatial refinement steps and  \abrevx{5} parametric refinement steps 
 (\abrevx{25} iterations in total).

\begin{figure}[!t]
\centering
\includegraphics[width = 0.45\textwidth]{{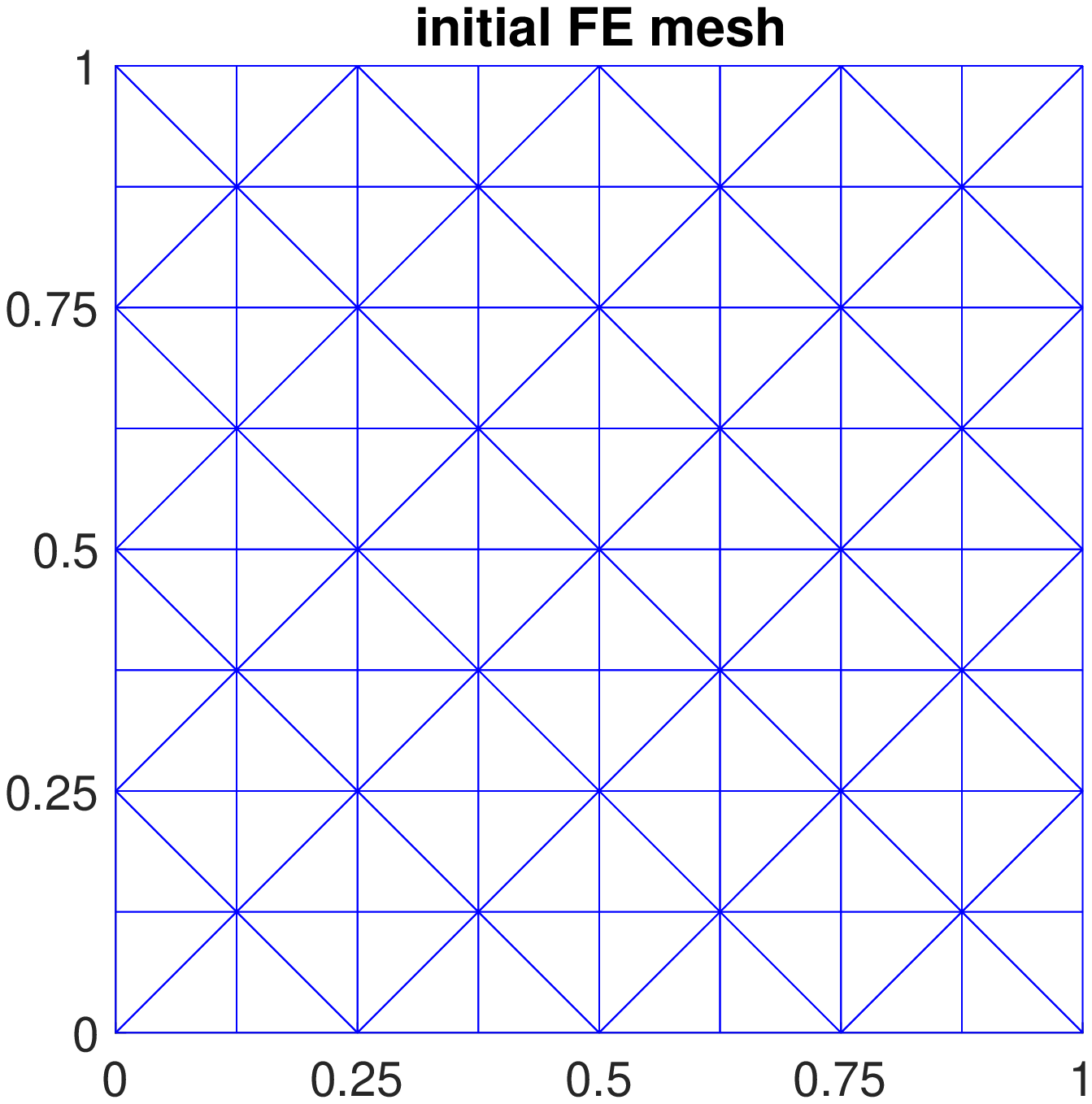}}
\includegraphics[width = 0.45\textwidth]{{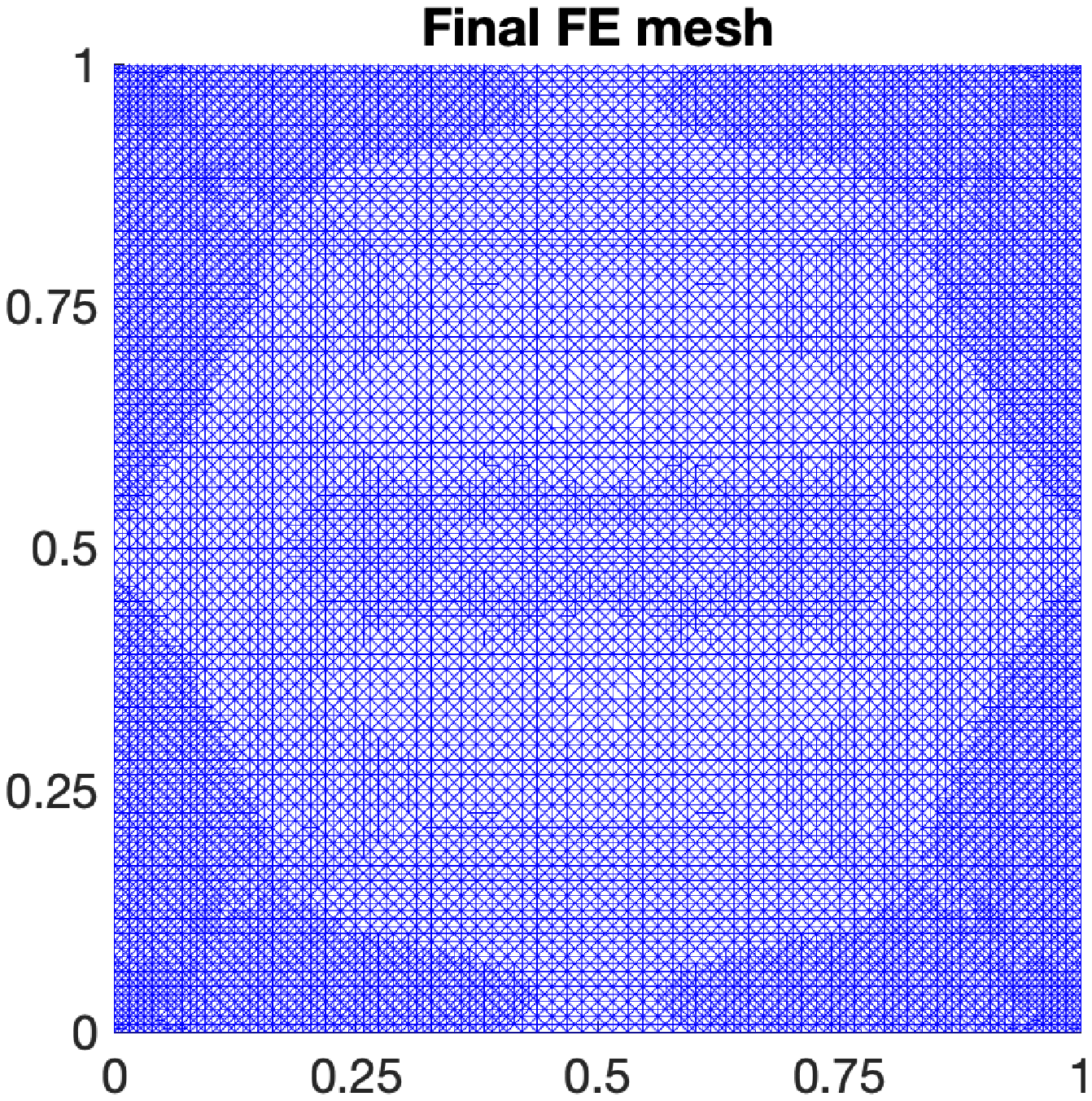}}
\caption{The initial mesh  (left) and  final mesh (right) generated by the 
single-level adaptive \rev{SC-FEM} strategy for test case I with error tolerance set to \abrevx{\tt 6e-3}.}
\label{fig:sc5.1mesh}
\end{figure}

The plots in Fig.~\ref{fig:sc5.1mesh} show the
initial mesh and the mesh  when the error tolerance is reached.
The final mesh can be seen to be  locally refined to resolve  weak
singularities in the corners. The number of vertices in the final mesh (that is,
the dimension of the linear system that is solved at every collocation point)
is \abrevx{\tt 16,473}.

\begin{table}[!t]
{\small
\begin{center}
\caption{Anisotropic approximation comparison  for test case I.}
\begin{tabular}{ccc|cc}
{\Large\strut}
 && sparse grid indices & parametric approximation space for SG & \\[3pt]
&&   1     1     1     1      &  0    0    0    0    0  & \\ 
&&   2     1     1     1     &   1    0    0    0    0  & \\ 
&&   3     1     1     1     &  2    0    0    0    0  & \\
 &&   4     1     1     1     &    3    0    0    0    0  & \\
 &&   5     1     1     1     &    4    0    0    0    0 & \\
 &&   1     2     1     1     &   0    1    0    0    0 & \\
 &&   2     2     1     1     &   1    1    0    0    0 & \\
 &&   3     2     1     1     &   2    1    0    0    0  & \\
 &&   1     3     1     1     &   0    2    0    0    0 & \\
 &&   2     3     1     1     &   1    0    0    1    0  & \\
 &&   3     3     1     1    &    1    0    1    0    0 & \\
&&    1     1     2     1    &    0    0    1    0    0 & \\
&&    1     1     3     1     &   2    0    1    0    0  & \\
&&                                &   3    1    0    0    0  & \\
&&                                &    0    0    0    1    0   & \\
&&                                &    0    0    0    0    1 & 
\end{tabular}
\end{center}
}
 \label{tab:indices}
\end{table}

The  parametric approximation initially consists of a single collocation point.  There 
were  \abrevx{13} Clenshaw--Curtis sparse grid collocation points when the error tolerance is reached.
The corresponding highly anisotropic sparse grid indices are listed in Table~1.
We note that  the resolution is concentrated in the first three coordinates and
that there is no refinement \abrevx{in the coordinates} corresponding to the  \abrevx{fourth mode}
in the expansion \eqref{kl}. (The  \abrevx{fourth mode} would be  activated  if we were to run with
\abrevx{a slightly tighter} tolerance,  say \abrevx{\tt 3e-3}.)

\begin{figure}[!th]
\centering
\includegraphics[width = 0.75\textwidth]{{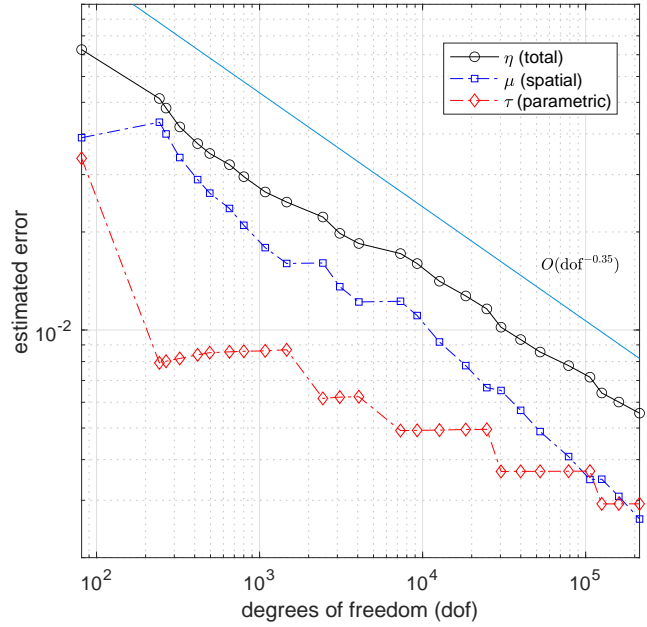}}
\includegraphics[width = 0.75\textwidth]{{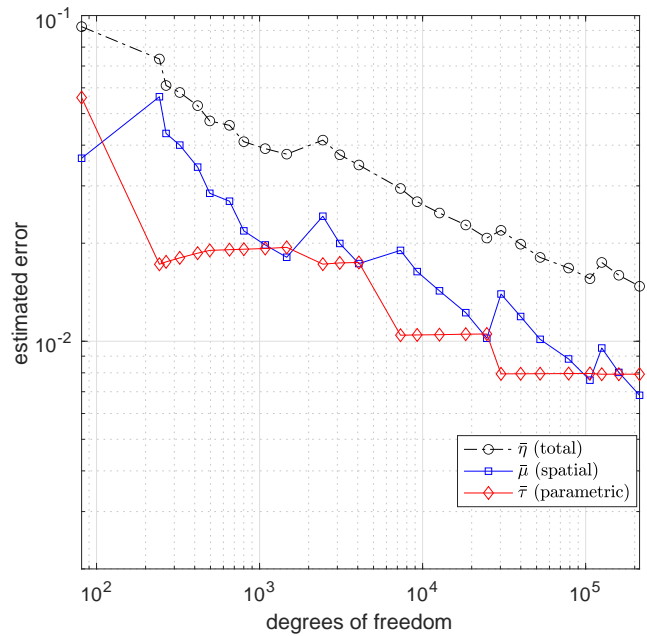}}
\caption{Evolution of the error estimates  (top) and  the error reduction indicators  (bottom) 
generated by the  single-level \rev{SC-FEM} adaptive strategy  for test case I
with error tolerance set to  \abrevx{\tt 6e-3}.}
\label{fig:sc5.1errors}
\end{figure}

The evolution of the component error estimates and the \abrevx{global error indicators}
is reported in  Fig.~\ref{fig:sc5.1errors}.  The key point here is that the nature  of  the refinement 
step (parametric or spatial)  is  determined by the relative size of the component global 
\abrevx{\it error indicators} (shown in the bottom plot). Thus it is reassuring to see 
the associated  parametric and spatial error estimates
(top plot)  decrease monotonically after the first few steps. Note that 
 if the saturation assumption \eqref{eq:saturation} is uniformly satisfied then the 
combined error estimate  \rev{$\eta_\ell := \mu_\ell  + \tau_\ell$} is guaranteed to 
decrease at every step.
In contrast, the total error  indicator \rev{$\bar\eta_{\ell} := \bar\mu_{\ell}  + \bar\tau_{\ell}$}
can be seen to {\it increase} at  iteration steps that follow a parametric refinement.  
The fact that the \abrevx{rate of the estimated error reduction in the bottom plot is 
 much slower \abrevx{than} the  rate of error reduction in the top plot} clearly  shows the  necessity 
 of computing \rev{$\eta_\ell$} periodically in Algorithm~\ref{algorithm}.

To check the robustness of the SC error estimation strategy we can compare the 
pattern of refinement  with the pattern that results when the same test problem \rev{(with $M = \infty$)} is 
solved   using the single-level stochastic Galerkin adaptive strategy \abrevx{in~\cite{br18, bprr18}}
that is built into T-IFISS~\cite{brs20}(with \abrevx{a slightly smaller accuracy tolerance}).   
When we ran this test,  the numerical solution generated by SG
is  visually identical to that reference solution  in Fig.~\ref{fig:sc5.1sol} with agreement to  \abrevx{4} decimal 
digits in the maximum value of the mean (\abrevx{\tt 0.07582} vs  \abrevx{\tt 0.07581})  \abrevx{as well as} 
 the maximum value of the standard deviation ({\tt 0.00710} vs  \abrevx{\tt 0.00709}).

\begin{figure}[!th]
\centering
\includegraphics[width = 0.75\textwidth]{{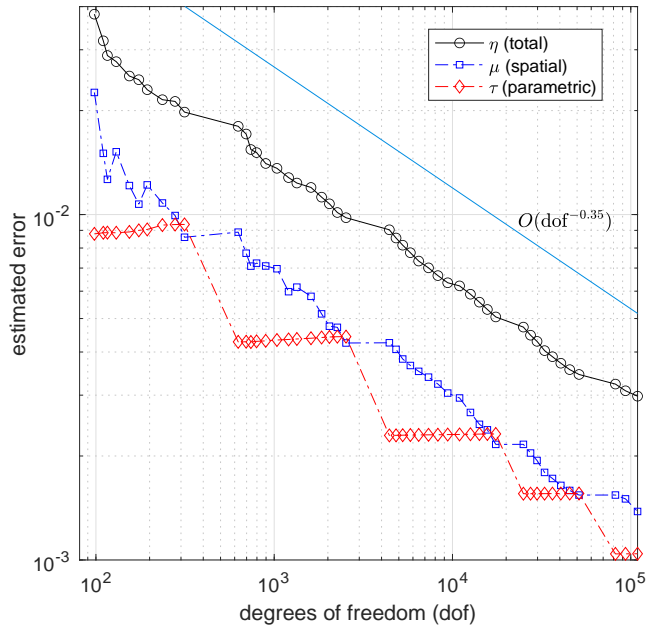}}
\caption{Evolution of the error estimates  and  the error reduction indicators 
generated by running single-level SG  for test case I with error tolerance 
set to {\tt 3e-3}.}
\label{fig:sg5.1errors}
\end{figure}
The evolution of the component SG error indicators and the global error estimate
is reported in  Fig.~\ref{fig:sg5.1errors}.  The total number of iterations is exactly twice
the number of adaptive SC steps\abrevx{,} with 4 parametric enrichment steps. 
The \abrevx{parametric approximation}
space \abrevx{obtained by SG} when the error tolerance \abrevx{is satisfied}  is listed in Table~1 
for comparison with SC. The \abrevx{main} difference is that there
are 5 parameters included in the SG approximation space  
when the algorithm terminated.\footnote{The
fourth parameter was  activated at step {\tt 37} and the fifth at step {\tt 45}.}
Comparing with adaptive SC the total number of degrees of freedom
was reduced by a factor of about 2 ({\tt 109,152} vs \abrevx{\tt 214,149 }) as was the
overall computation time (29 seconds vs 71 seconds). Reassuringly, the rates
of convergence of the SC and SG algorithms can be seen to be \abrevx{closely
matched.}

\begin{figure}[!t]
\centering
\includegraphics[width = 0.45\textwidth]{{sc5.1mesh0x.eps}}
\includegraphics[width = 0.45\textwidth]{{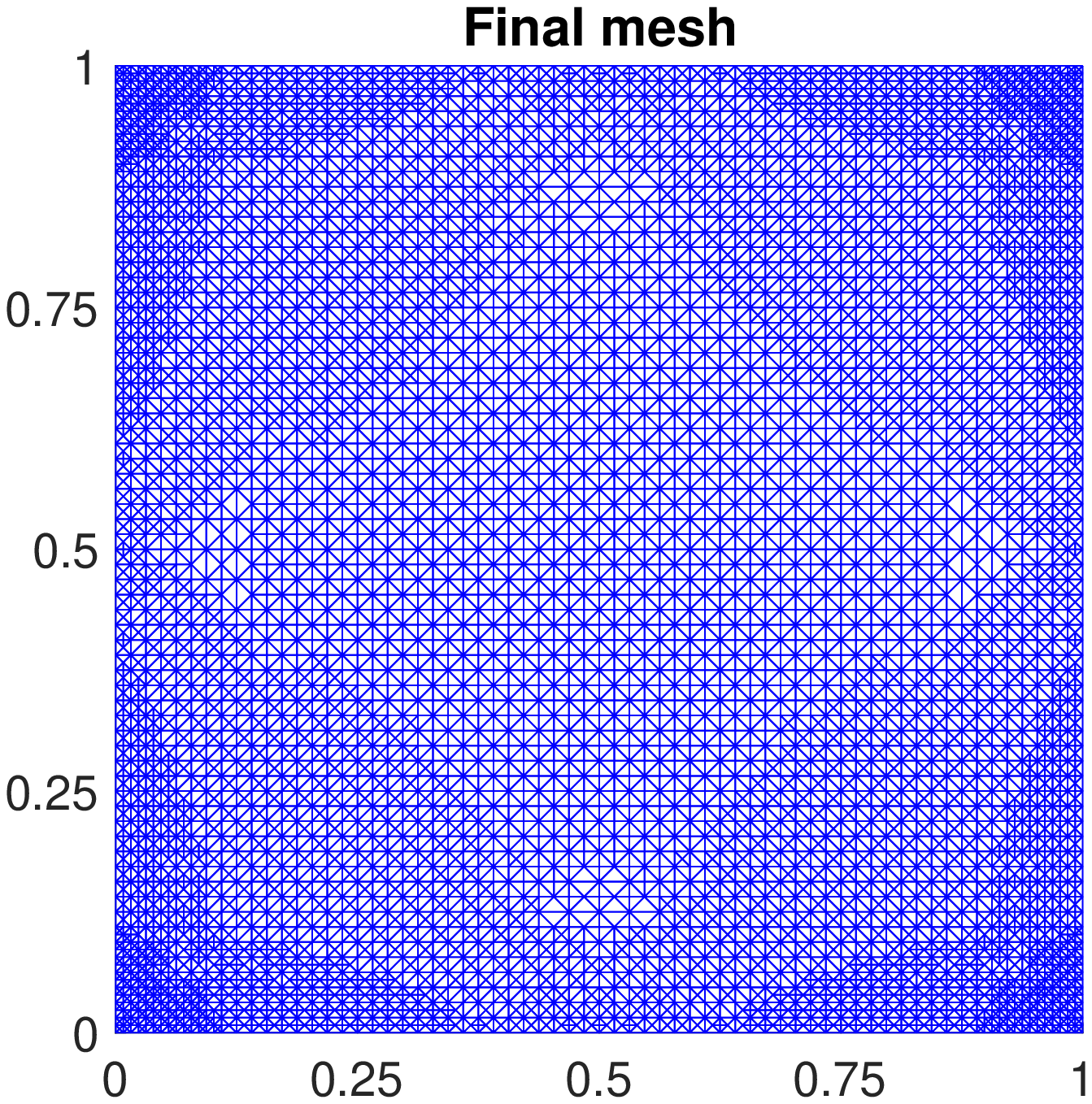}}
\caption{The initial mesh  (left) and  final mesh (right) generated by  running
single-level adaptive SG  for test case I with error tolerance set to {\tt 3e-3}.}
\label{fig:sg5.1mesh}
\end{figure}

The plots in Fig.~\ref{fig:sg5.1mesh} show the
initial mesh together with the mesh  when the SG error tolerance is reached.
 The number of vertices in the final mesh is {\tt 7,134}. This is less
 refined than the final mesh in Fig.~\ref{fig:sc5.1mesh} as might be expected.
 One obvious  difference between the two final meshes is the fact that the mesh 
 generated by adaptive SC has local resolution that captures the 
 distribution of  the variance in the reference solution (cf. Fig.~\ref{fig:sc5.1sol}).

\subsection{Test case II:  nonaffine coefficient data}\label{sec:nonaffineresults}

In this case, we set $f = 1$ and solve  the model problem
on the L-shaped domain \rev{$D = (-1, 1)^2\backslash (-1, 0]^2$} with
coefficient $a(x, \y) = \exp(h( x, \y))$,
where the exponent field  $h(x, \y)$ has affine dependence on
parameters $y_m$  that are images of uniformly 
distributed independent mean-zero random variables,
\begin{align} \label{kll}
h(x, \y) = h_0(x) + \sum_{m = 1}^M h_\abrevx{m}(x) \, y_m,\quad
   x \in D,\ \y \in \Gamma.
\end{align}
We further specify $h_0(x) \,{=}\, 1$ and 
$h_m(x) \,{=}\, \sqrt{\lambda_\abrevx{m}} \varphi_m(x)$
($m \,{=}\, 1,\ldots,M$). Here $\{(\lambda_m, \varphi_m)\}_{m=1}^\infty$ 
are the eigenpairs of the integral operator
$\int_{\rev{ D \cup (-1,0]^2} } \hbox{\rm Cov}[\rev{h}](x, x') \varphi(x')\, \hbox{d} x' $
with a  synthetic covariance function given by
\begin{align} \label{cov}
\hbox{\rm Cov}[\rev{h}](x, x') 
= \sigma^2 \exp
\left( -\frac{| x_1 - x_1' |}{\ell_1}  -\frac{| x_2 - x_2' | }{\ell_2} \right),
\end{align}
where $\sigma$ is the standard deviation 
and $\ell_1$, $\ell_2$ are correlation lengths
(we set $\ell_1 = \ell_2 = 1$). The resulting parametric model problem is uniformly  
well posed in the sense that \eqref{eq:amin:amax} is satisfied for any choice of $M$.
If $M$ is fixed then  the challenge is to retain  
robustness when  the standard deviation is increased. A test case  such as this is 
more amenable to SC-FEM approximation than the affine test case discussed above. 
The sparsity of the linear algebra would be severely compromised if this test problem were
solved using a  stochastic Galerkin approximation strategy\abrevx{; see, e.g.,~\cite{bx2020}}.

We present results for three test problems associated with different combinations of the 
number of parameters $M$ and the standard deviation $\sigma$: (a) $M = 4$ 
and $\sigma = 0.5$, (b) $M = 8$ $\sigma = 0.5$ and (c) $M = 4$ and $\sigma = 1.5$. 
For all tests we specify the same tolerance (\abrevx{\tt 6e-3}) and 
run the algorithm with marking parameters $\theta_{\rev{\X}} = \theta_{\rev{\Colpts}} = 0.3$.
To assess the effectivity of the error estimation process
we also computed a reference solution $u_{\text{ref}}$ 
as a proxy of the exact solution to each problem.
The reference solution(s) were generated using the minimum isotropic index set 
containing the final index set  from the adaptive computation
together with a piecewise {\it quadratic} finite element approximation space 
defined on the final mesh from the adaptive computation.

To measure the quality of the error estimate \abrevx{$\eta_\ell = \mu_\ell + \tau_\ell$} 
we compute  an effectivity index at each iteration via
\begin{align}
\Theta_\ell  = \abrevx{ \frac{\eta_\ell}{\| u_{\text{ref}} \,{-}\, \rev{u_{\ell}^{\rm SC}} \|} }.
\end{align}
Thus, $\Theta_\ell$ being close to $1$ suggests that \abrevx{$\eta_\ell$}  is an effective
estimate of the norm of the error.

 \begin{figure}[!th]
\centering
\includegraphics[width = 0.75\textwidth]{{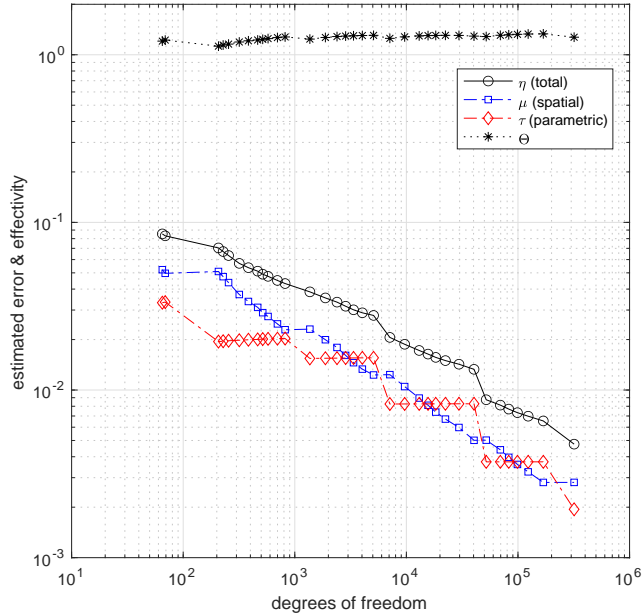}}
\caption{Evolution of the error estimates  generated by running 
the single-level adaptive strategy for test case II with $M=4$ and  $\sigma = 0.5$
 with the error tolerance  set to \abrevx{\tt 6e-3}.}
\label{fig:sg5.2aerrors}
\end{figure}

\begin{figure}[!t]
\centering
\includegraphics[width = 0.45\textwidth]{{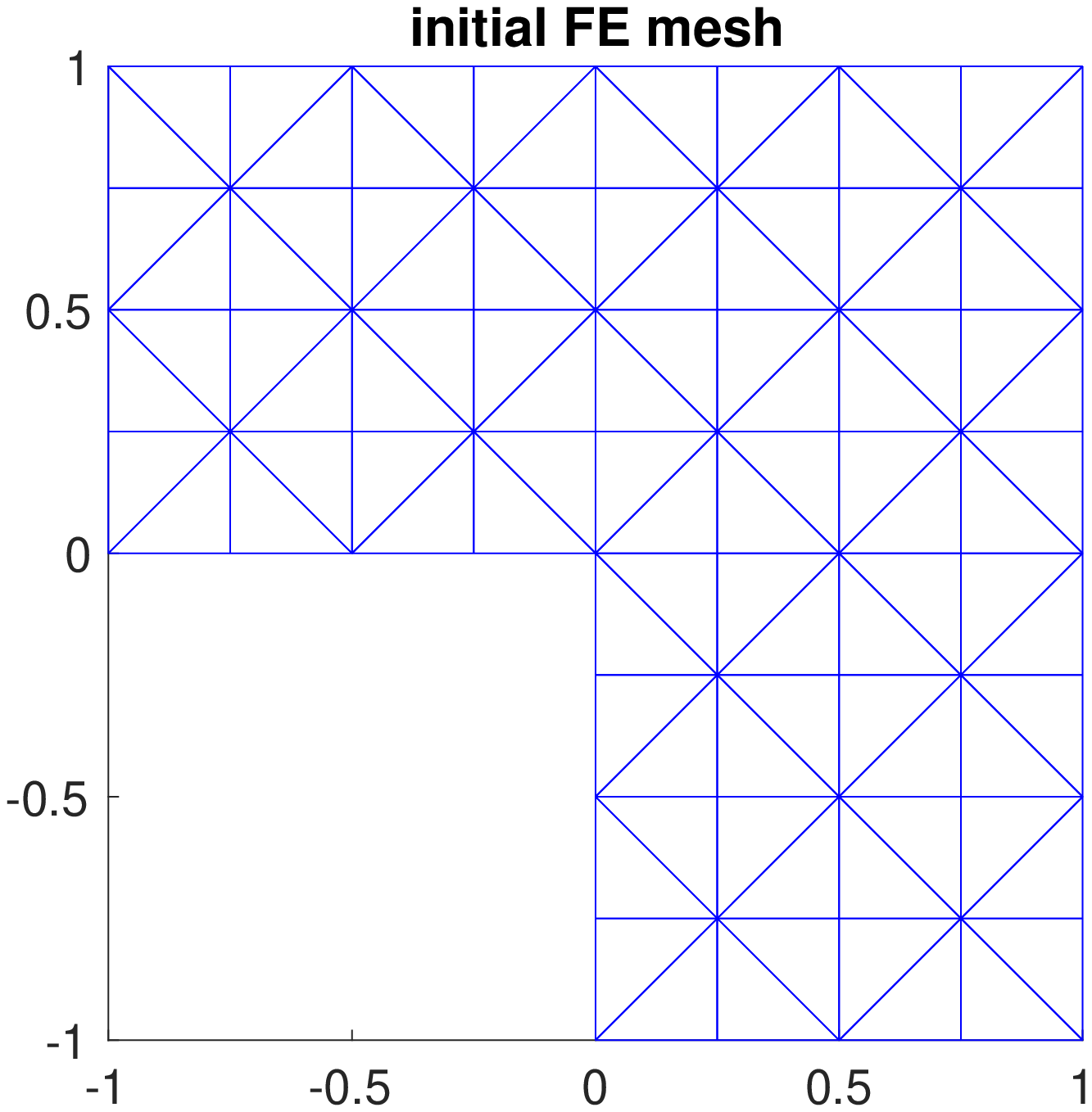}}
\includegraphics[width = 0.45\textwidth]{{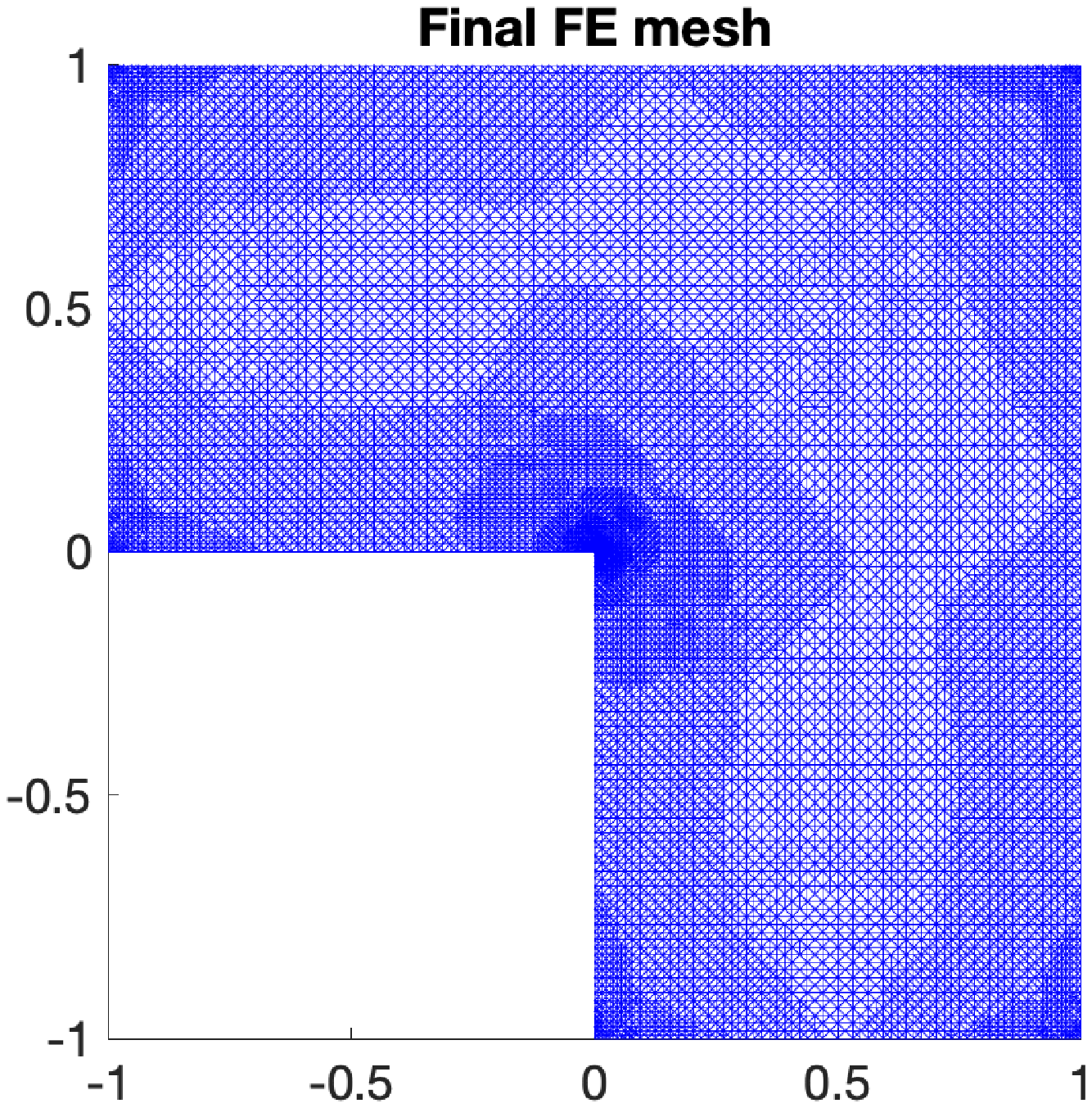}}
\caption{The initial mesh  (left) and  final mesh (right) generated by  running
the single-level adaptive strategy for test case II  with $M=4$ and  $\sigma = 0.5$
 with error tolerance set to \abrevx{\tt 6e-3}.}
\label{fig:sc5.2mesh}
\end{figure}

We record the evolution of the error estimates and the component 
spatial and parametric \rev{contributions} together with the effectivity index at each 
iteration of the single-level adaptive algorithm. Results for the 
 first parameter combination are shown in Fig.~\ref{fig:sg5.2aerrors}.  The error tolerance 
 was satisfied after 33 iterations which included  5 parametric refinement steps. 
 The plots in Fig.~\ref{fig:sc5.2mesh} show the
initial mesh and the mesh  when the error tolerance is reached.
The degree of local refinement in  the final mesh is strongest around
the reentrant corner  but is noticeable  at all  corners of the domain.
The number of vertices in the final mesh is {\tt 18,737}\abrevx{,}
and \abrevx{there are} 17 collocation points in the final \abrevx{sparse grid}.
The computed effectivity indices \abrevx{plotted} in  Fig.~\ref{fig:sg5.2aerrors} can
be seen to stay close to unity throughout, ranging from a minimum value  of 1.125 to a 
maximum value of 1.333.

 \begin{figure}[!th]
\centering
\includegraphics[width = 0.75\textwidth]{{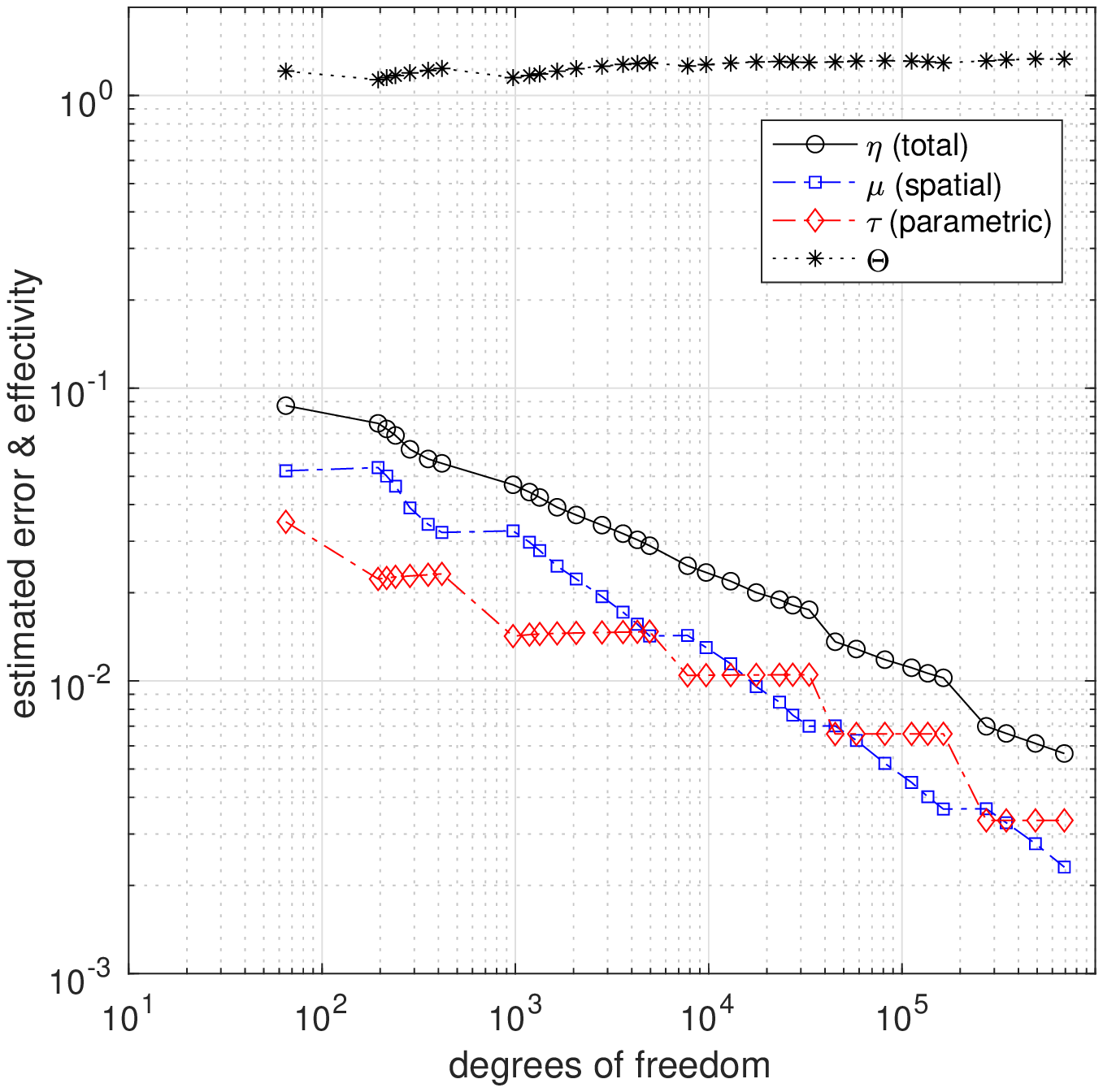}}
\caption{Evolution of the error estimates  generated by running 
the single-level adaptive strategy for test case II with $M=8$ and  $\sigma = 0.5$
 with the error tolerance  set to \abrevx{\tt 6e-3}.}
\label{fig:sg5.2berrors}
\end{figure}

\begin{figure}[!th]
\centering
\includegraphics[width = 0.75\textwidth]{{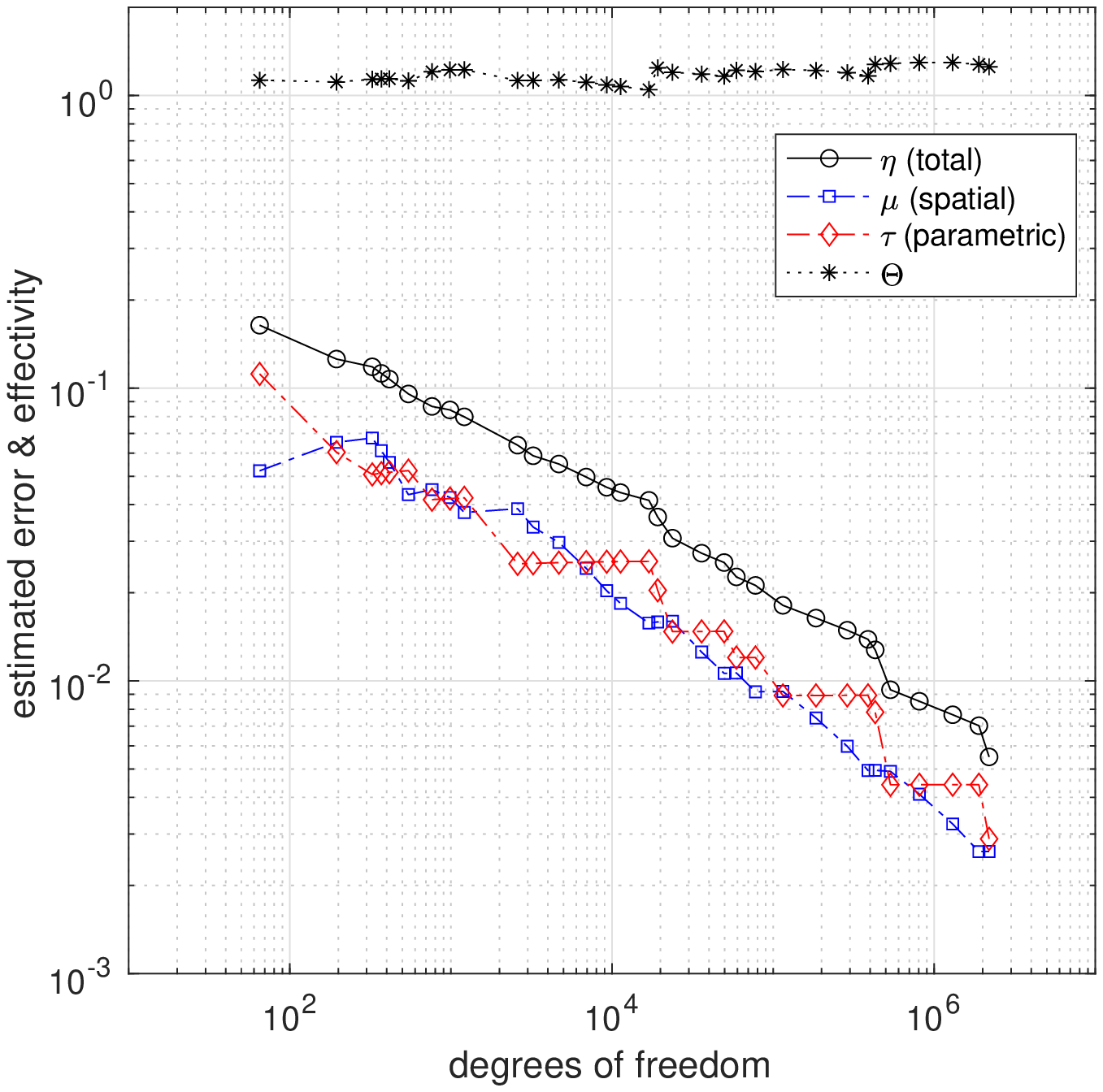}}
\caption{Evolution of the error estimates  generated by running 
the single-level adaptive strategy for test case II with $M=4$ and  $\sigma = 1.5$
 with the error tolerance  set to \abrevx{\tt 6e-3}.}
\label{fig:sg5.2cerrors}
\end{figure}

Convergence histories for \abrevx{test problems with the other two parameter combinations} are
presented in Fig.~\ref{fig:sg5.2berrors} and in Fig.~\ref{fig:sg5.2cerrors}.
While the number of iteration steps  was reduced from 33 to 32 when $M$ was
increased from 4 to 8, the number of parametric refinement steps did not change.
The number of degrees of freedom at the end is noticeably larger however: this is mainly
because  the number of vertices in the final mesh has increased 
 by 50\% (from  {\tt 18,737} to {\tt 27,651}).
Reassuringly, the number of  collocation points in the 
final \abrevx{sparse grid} has been kept under control, increasing from 15  in $\R^4$ to 25 in $\R^8$.
(The reason for this is  that  the original $M=4$ sparse grid index set has simply been augmented by 
extra index entries generating \rev{nodes} along the \rev{axes}  of the  4 additional coordinate directions.
The \rev{nodes} in the original directions were unchanged.)

The key point that is worth reiterating  is the following: while
the \rev{total error indicator} \rev{$\bar\eta_{\ell} = \bar\mu_{\ell}  + \bar\tau_{\ell}$}
is not  robust  as a measure of the discretization error in cases like this where
 $M$ is increased, the relative size of the components still seem to provide reliability in 
the refinement pattern. Moreover, looking at the results in Fig.~\ref{fig:sg5.2berrors}
the effectivity indices can be seen to stay close to unity throughout the adaptive process,
ranging from a minimum value  of 1.132 to a maximum value of 1.332.  Thus our 
expectation is  that a stopping criterion  based on \rev{$\eta_\ell$} will be reliable  
in general, not just when $M$ is small. 

Turning to the \abrevx{test problem for the third parameter combination}, wherein the 
standard deviation is increased while keeping $M$ fixed, we again see a significant  increase in 
the dimension of the discrete problem at the termination of the adaptive process. While the
number of iteration steps is reduced from 33 to 31, the pattern
of the refinement is very different in this case. Looking carefully at
the convergence history reproduced  in Fig.~\ref{fig:sg5.2cerrors}\abrevx{,}  a total of 11 parametric 
refinement steps  can be identified. As a result,  the number  of   collocation points 
in the final sparse grid  index set  is  significantly increased (mirroring the
increase in the uncertainty) from 15  when $\sigma$ is  0.5 
to 59  when $\sigma$ is increased to 1.5.
Thus, since the number of vertices in the final mesh has also doubled
 (going from  {\tt 18,737} to {\tt 37,133}) there is an order of magnitude
 increase in  the total number of degrees of freedom needed to solve the \abrevx{test}
 problem to the specified accuracy.
 \abrevx{The effectivity of the error estimation is also retained, with  indices 
  staying between 1.047 and 1.296.}

\abrevx{
The numerical results presented above demonstrate the effectivity and the robustness
of our distinctive SC error estimation strategy as well as the utility of the error indicators guiding 
the adaptive refinement process. Optimality of convergence is, however, precluded when using the 
single-level approach. In Part II of this work, we will investigate the realization of close-to-optimal 
convergence rates within a multilevel framework.
}

\bibliographystyle{siam}
\bibliography{references}

\end{document}